\RequirePackage{fix-cm}
\documentclass[smallextended]{svjour3}       % onecolumn (second format)
\smartqed  % flush right qed marks, e.g. at end of proof
\usepackage{amssymb}  
\usepackage{float} %%float: \begin{figure}[H]
\usepackage{amsmath}
\usepackage{graphicx}
\usepackage{subfig}
\usepackage{color}
\begin{document}

\title{Some algorithms for the mean curvature flow under topological changes%\thanks{Grants or other notes
%about the article that should go on the front page should be
%placed here. General acknowledgments should be placed at the end of the article.}
}
%\subtitle{Do you have a subtitle?\\ If so, write it here}

%\titlerunning{Short form of title}        % if too long for running head

\author{Arthur Bousquet         \and
        Yukun Li  \and  
        Guanqian Wang%etc.
}

%\authorrunning{Short form of author list} % if too long for running head

\institute{A. Bousquet \at
              Department of Mathematics, Lake Forest College \\
              %Tel.: +123-45-678910\\
              %Fax: +123-45-678910\\
              \email{bousquet@mx.lakeforest.edu}           %  \\
%             \emph{Present address:} of F. Author  %  if needed
           \and
           Y. Li \at
              Department of Mathematics, University of Central Florida, Orlando, Florida, United States\\
              \email{yukun.li@ucf.edu}  
                         \and
           G. Wang \at
              Department of Mathematics, University of Central Florida, Orlando, Florida, United States\\
              \email{guanqian.wang@ucf.edu}  
}

\date{Received: date / Accepted: date}
% The correct dates will be entered by the editor

\maketitle

\begin{abstract}
This paper considers and proposes some algorithms to compute the mean curvature flow under topological changes. Instead of solving the fully nonlinear partial differential equations based on the level set approach, we propose some minimization algorithms based on the phase field approach. It is well known that zero-level set of the Allen-Cahn equation approaches the mean curvature flow before the onset of the topological changes; however, there are few papers systematically studying the evolution of the mean curvature flow under the topological changes. There are three main contributions of this paper. First, in order to consider various random initial conditions, we design several benchmark problems with topological changes, and we find different patterns of the evolutions of the solutions can be obtained if the interaction length (width of the interface) is slightly changed, which is different from the problems without topological changes. Second, we propose an energy penalized minimization algorithm which works very well for these benchmark problems, and thus furthermore, for the problems with random initial conditions. Third, we propose a multilevel minimization algorithm. This algorithm is shown to be more tolerant of the unsatisfying initial guess when there are and there are no topological changes in the evolutions of the solutions.
\keywords{mean curvature flow \and Allen-Cahn equation \and topological changes
\and energy penalized minimization algorithm \and multilevel minimization algorithm}
% \PACS{PACS code1 \and PACS code2 \and more}
\subclass{65N12 \and 65N22 \and 65N30 \and 65N55}
\end{abstract}

\section{Introduction}\label{sec1}
The moving interface problems refer to a type of problems with interfaces which move in point-wise velocity. Among various kinds of moving interface problems, the mean curvature flow is the most basic and important one. It is defined to be the case when the outward normal velocity is equal to the negative point-wise mean curvature, i.e.,
\begin{equation}\label{eq20190804_1}
V_n(x,t)=-H(x,t)\qquad x\in\Gamma_t,
\end{equation}
where $V_n(x,t)$ and $H(x,t)$ denote the outward normal velocity and the mean curvature of the interface $\Gamma_t$ respectively. The moving interface problems have been broadly applied to a few different areas such as fluid mechanics, gas dynamics, biology, financial mathematics, etc. There are basically two approaches to formulate the moving interface problems. The first approach is the direct approach, i.e., it tracks the interface directly, including the parametric finite element method \cite{barrett2007parametric}, the volume of fluid method \cite{hirt1981volume}, the immersed boundary method \cite{peskin2002immersed}, the front tracking method \cite{unverdi1992front}, the immersed interface method \cite{leveque1994immersed}, and so on. The second approach is the indirect approach, i.e., it seeks another quantity instead of tracking the interface directly, including the level set method \cite{osher1988fronts} and the phase field method \cite{rayleigh1892xx}. The biggest advantage of the second approach is {\color{black} that} it can easily handle the topological changes compared to the first approach. 

The objective of the second approach is to seek an indirect quantity called the level set function (the level set method) or the phase field function (the phase field method). The level set formulation of the mean curvature flow \eqref{eq20190804_1} can be written as \cite{chen1991uniqueness,evans1992phase,evans1992motion}
\begin{align}
\omega_t&=\Delta\omega-\frac{D^2\omega\nabla\omega\cdot\nabla\omega}{|\nabla\omega|^2},\label{eq20190804_2}\\
\omega(\cdot,0)&=\omega_0(\cdot),\label{eq20190804_3}
\end{align}
where $D^2\omega$ denotes the Hessian of $\omega$ and $\omega_0$ satisfies $\Gamma_0=\{x\in\mathbb{R}^d;\ \omega_0(x)=0\} (d=2,3)$. {\color{black}This means the interface $\Gamma_t$ evolved by the mean curvature flow is the zero-level set of the solution $\omega$.} It is explained in \cite{evans1992motion} that problem \eqref{eq20190804_2}--\eqref{eq20190804_3} has a unique continuous solution and the level set of $\omega$ evolves according to the mean curvature flow. Some numerical methods are proposed to discretize \eqref{eq20190804_2}--\eqref{eq20190804_3} in \cite{deckelnick2003numerical,m2017approximations,kovacs2018convergent,walkington1996algorithms}. The phase field formulation of the mean curvature flow \eqref{eq20190804_1} is the singular perturbation of the heat equation called the Allen-Cahn equation with the following initial and boundary conditions
\begin{alignat}{2}
u_t-\Delta u+\frac{1}{\epsilon^2}f(u)&=0  &&\qquad \mbox{in } 
\Omega_T:=\Omega\times(0,T),\label{eq1.1}\\
\frac{\partial u}{\partial n} &=0 &&\qquad \mbox{in }
\partial\Omega_T:=\partial\Omega\times(0,T), \label{eq1.2}\\
u &=u_0 &&\qquad \mbox{in }\Omega\times\{t=0\}. \label{eq1.3}
\end{alignat}
where $\epsilon$ denotes the interaction length, $\Omega\subseteq \mathbb{R}^d$ is a bounded domain, $f=F'$ for some double well potential function $F$, and $F(u)=\frac{1}{4}(u^2-1)^2$ is used in this paper. The zero-level sets of $u$ are proved to converge to the mean curvature flow in \cite{evans1992phase,ilmanen1993convergence}, and the zero-level sets of numerical solutions are proved to converge to the mean curvature flow in \cite{Feng_Prohl03,Feng_Wu05,feng2014analysis,kessler2004posteriori,li2015numerical,xu2019stability}. We refer the readers to further numerical discussions on the stochastic mean curvature flow \cite{feng2014finite,feng2017finite,feng2018strong,majee2018optimal} and the Hele-Shaw flow \cite{feng2016analysis,li2017error,wu2018analysis} {\color{black}which is an another} fundamental moving interface problem.
 
The above results discuss the evolutions of the solutions without topological changes. In \cite{bartels2011quasi,bartels2011robust,bartels2010lower}, the authors consider the topological changes under some initial conditions, but the evolutions of the solutions behave similarly under different values of the interaction length. In \cite{du2005retrieving}, the authors suggest using Euler number to retrieve topological information thus to avoid unphysical changes of the topology. Different from the above results, this paper considers some sensitive benchmark initial conditions, i.e. the topology or the patterns of the evolutions of the solutions could be totally different under different values of the interaction length. The motivation of investigating these sensitive initial conditions is that they frequently happen when the initial conditions are random or rough. In this paper, we assume there is no fattening phenomenon \cite{cesaroni2018fattening,evans1992motion,ilmanen1993convergence} in the evolutions of the solutions under specially designed initial conditions. The objectives of this paper are threefold. First, the random initial conditions can be used to test the topological changes in the worst scenarios, so this paper constructs and solves several these kinds of sensitive benchmark problems, and finally tests the random initial conditions. Second, this paper proposes an energy penalized minimization algorithm which is aimed to accurately track the mean curvature flow under topological changes and efficiently solve the problem, i.e., it allows big time step size compared to classical numerical methods for the phase field models. All the proposed benchmark problems as well as the problems with random initial conditions will be used to validate this algorithm. Third, a multilevel minimization algorithm is proposed to handle bad initial guesses and the topological changes. The objective of this algorithm is to construct convex functionals such that their minimization problems {\color{black}do not significantly} depend on the initial guesses.

This paper has four additional sections. In section 2, we construct several benchmark problems with topological changes. The level set method, the phase field method, and the energy minimization method (without penalty) are used to compare the evolutions of the solutions. In section 3, we propose the energy penalized minimization algorithm. This algorithm is implemented for the aforementioned benchmark problems as well as problems with random initial conditions. In section 4, we propose a multilevel minimization algorithm. Some examples with and without topological changes are used to justify this algorithm. In section 5, we give a short concluding remark to summarize the results in this paper. 

\section{Some Benchmark Problems.}\label{sec2}
In this section, we will state a few benchmark problems of the mean curvature flow under topological changes. The level set method, the energy minimization method, and the phase field method are employed in this section. The solution computed by the level set method is considered the reference solution. We begin with the discretization of the phase field model. We first introduce some notations. Let $\mathcal{T}_h$ be a
shape-regular triangulation of $\Omega \subset \mathbb{R}^d \ (d=2,3)$, $K\in \mathcal{T}_h$ represent each element, and $h$ denote the diameter of $K$. Define the finite element
space $V_h$ by
\begin{equation} \label{fem_space}
V_h = \bigl\{v_h\in C(\overline{\Omega}):\ v_h|_K\in P_r(K)\bigr\},
\end{equation}
where $C(\overline{\Omega})$ denotes the set of all continuous functions on $\overline{\Omega}$ and $P_r(K)$ denotes the set of all polynomials whose degrees are less than or equal to a given positive integer $r$ on $K$.  {\color{black}Let $k$ be the time step size,} $\|\cdot\|_{L^2}$ be the $L^2$-norm, and $(\cdot, \cdot)$ be the $L^2$-inner product over the domain $\Omega$. {\color{black}The uniform time step size $k$ can be extended to non-uniform time step size $k_n$, $n=1,2,\cdots
$ in this paper.} The first well-known scheme is the standard first-order fully implicit scheme (FIS), which is defined by seeking $u_h^n\in V_h$ for $n=1,2,\cdots
$, such that
\begin{equation}\label{FIS-AC}
(\frac{u_h^{n}-u_h^{n-1}}{k},v_h)+(\nabla u_h^{n},\nabla v_h) +
\frac{1}{\epsilon^2}((u_h^{n})^3-u_h^n,v_h) = 0 \qquad\forall v_h\in V_h.
\end{equation}
The following are some usual numerical schemes which will be used in this paper. The convex splitting scheme is to find $u_h^n\in V_h$ for $n=1,2,\cdots$, such that
\begin{equation}\label{css}
(\frac{u_h^{n}-u_h^{n-1}}{k},v_h)+(\nabla u_h^{n},\nabla
v_h)+\frac{1}{\epsilon^2}((u_h^{n})^3-u_h^{n-1},v_h)=0
\qquad\forall v_h\in V_h.
\end{equation}
The Semi-implicit scheme is to find $u_h^n\in V_h$ for $n=1,2,\cdots$, such that
\begin{align}\label{semi}
(\frac{u_h^{n}-u_h^{n-1}}{k},v_h)&+(\nabla u_h^{n},\nabla
v_h)\\
&+\frac{1}{\epsilon^2}((u_h^{n-1})^3-u_h^{n-1},v_h)=0
\qquad\forall v_h\in V_h.\notag
\end{align}
The modified Crank-Nicolson scheme is to find $u_h^n\in V_h$
for $n = 1, 2, \cdots$, such that
\begin{align}\label{eq:CN-AC}
\bigl(\frac{ u_h^{n}- u_h^{n-1}}{k},v_h\bigr)&+\bigl(\frac{\nabla
    u_h^{n} + \nabla u_h^{n-1}}{2}, \nabla v_h\bigr)\\
&+\frac{1}{\epsilon^2}(\tilde F[u_h^n, u_h^{n-1}], v_h) = 0 \quad \forall
v_h\in V_h,\notag
\end{align}
where
$$
\tilde F[u, u_h^{n-1}]=
\begin{cases}
\frac{F(u)-F(u_h^{n-1})}{ u - u_h^{n-1}} & u\neq u^{n-1}_h,\\
u^3 - u & u =  u^{n-1}_h.
\end{cases}
$$
A review of different numerical schemes can also be found in \cite{shen2011spectral}, and a few benchmark problems can be found in \cite{church2019high}.

Next we define the following free-energy functional $J_\epsilon^{\rm AC}$ and the  discrete energy $E_n^{\rm AC}$
\begin{align}\label{eq20200423_2}
J_\epsilon^{\rm AC}(v) &:= 
\int_\Omega \Bigl( \frac12 |\nabla v|^2
+ \frac{1}{\epsilon^2} F(v) \Bigr)\, dx,\\
E_n^{\rm AC}(u_h ;u_h^{n-1}) 
&:=J_\epsilon^{\rm AC}(u_h)+\frac{1}{2k}\int_{\Omega}(u_h - u_h^{n-1})^2dx.
\label{AC_energy}\end{align}
It is known that the Allen-Cahn equation is the $L^2$-gradient flow of the functional $J_\epsilon^{\rm AC}$, and it is proved that the fully implicit scheme \eqref{FIS-AC} is equivalent to the following energy minimization problem:  
\begin{equation} \label{tra_energy}
u_h^n = \underset{u_h\in V_h}{\mathrm{argmin}} E_n^{\rm AC}(u_h;u_h^{n-1}).
\end{equation}
See \cite{xu2019stability} for the relations between other numerical schemes and the energy minimization methods.

When there are no topological changes, it was numerically observed in \cite{feng2014analysis,li2015numerical,zhang2009numerical} that the energy decays slowly. When there are topological changes, it was numerically observed that the energy decays fast. % and the spectrum of the Allen-Cahn operator blows up to negative infinity \cite{bartels2011quasi, bartels2011robust, bartels2010lower}.
It is very hard or even impossible to calculate the exact energy especially when the geometry is complex. Here we provide an energy equality to formally explain the above observations; {\color{black}i.e., when topological changes happen, the term $k\sum_{n=0}^\ell R^n_{\epsilon}$ in Theorem \ref{lem20200421} is larger since the profile of the solutions change more.} Note the discrete energy inequality of the Allen-Cahn equation was established in \cite{feng2014analysis,xu2019stability}.
\begin{theorem}\label{lem20200421}
Let $u_h^n$ be a solution of scheme \eqref{FIS-AC} {\color{black}and define the difference operator $d_t$ by 
\[d_t u^m:= \frac{u^m-u^{m-1}}{k}, m=1,2,\cdots,M,\] 
where $M$ is the largest positive integer less than or equal to $\frac{T}{k}$}. Then the following discrete energy equality holds:
\begin{equation*}
J_{\epsilon}^{AC}(u_h^\ell) + k\sum_{n=0}^\ell R^n_{\epsilon} = J_\epsilon^{AC}(u_h^0) 
\qquad\text{for } 0\leq \ell \leq M,
\end{equation*}
where the functional $J_\epsilon^{\rm AC}$ is defined in \eqref{eq20200423_2}, and $R^n_{\epsilon}$ is defined by
\begin{equation*}
R^n_{\epsilon} :=\Bigl(1- \frac{k}{2\epsilon^2} \Bigr)
\|d_t u_h^{n}\|_{L^2}^2 
+\frac{k}{2} \|\nabla d_t u_h^{n}\|_{L^2}^2
+\frac{k}{4\epsilon^2}\|d_t(u_h^{n})^2\|_{L^2}^2+\frac{k}{2\epsilon^2}\|u_h^{n}d_tu_h^{n}\|_{L^2}^2. 
\end{equation*}
\end{theorem}

\begin{proof}
Choosing $v_h=u_h^{n}-u_h^{n-1}$ in \eqref{FIS-AC}, we get
\begin{equation}\label{eq20200423_1}
\frac{1}{k}\|u_h^{n}-u_h^{n-1}\|_{L^2}^2+(\nabla u_h^{n},\nabla (u_h^{n}-u_h^{n-1})) +\frac{1}{\epsilon^2}((u_h^{n})^3-u_h^n,u_h^{n}-u_h^{n-1}) = 0.
\end{equation}
Define the backward difference operator   
\begin{align}\label{eq20200423_3}
d_t v^m:= \frac{v^m-v^{m-1}}{k} \qquad \forall m=1,2,\cdots,\ell.
\end{align}
Note the following identities
\begin{align}\label{eq20200423_4}
(\nabla u_h^{n},\nabla(u_h^{n}-u_h^{n-1}))=&\frac12\|\nabla u_h^{n}\|_{L^2}^2-\frac12\|\nabla u_h^{n-1}\|_{L^2}^2\\
&+\frac12\|\nabla (u_h^{n}- u_h^{n-1})\|_{L^2}^2,\notag\\
((u_h^{n})^3-u_h^n,u_h^{n}-u_h^{n-1})=&\frac{k}{4}d_t\|(u_h^{n})^2-1\|_{L^2}^2+\frac{k^2}{4}\|d_t(u_h^{n})^2\|_{L^2}^2\label{eq20200423_5}\\
&+\frac{k^2}{2}\|u_h^{n}d_tu_h^{n}\|_{L^2}^2-\frac{k^2}{2}\|d_tu_h^{n}\|_{L^2}^2.\notag
\end{align}
Then we have
\begin{align}\label{eq20200423_6}
&k\|d_tu_h^{n}\|_{L^2}^2+\frac12\|\nabla u_h^{n}\|_{L^2}^2-\frac12\|\nabla u_h^{n-1}\|_{L^2}^2+\frac12\|\nabla (u_h^{n}- u_h^{n-1})\|_{L^2}^2\\
&\qquad+\frac{k}{4\epsilon^2}d_t\|(u_h^{n})^2-1\|_{L^2}^2+\frac{k^2}{4\epsilon^2}\|d_t(u_h^{n})^2\|_{L^2}^2\notag\\
&\qquad+\frac{k^2}{2\epsilon^2}\|u_h^{n}d_tu_h^{n}\|_{L^2}^2-\frac{k^2}{2\epsilon^2}\|d_tu_h^{n}\|_{L^2}^2=0.\notag
\end{align}
Taking the summation on both sides of \eqref{eq20200423_6} over $n$ from $0$ to $\ell$ yields the conclusion.
\end{proof}

In the following tests, we compare the evolutions of the solutions of different numerical schemes, and the evolutions of the solutions of different methods such as the level set method, the phase field method and the energy minimization method. Specifically, in Test 1, we verify that all these popular numerical schemes behave similarly as long as the time step size is small enough. Therefore, we choose the fully implicit scheme in the following tests; in Test 2, we choose two circles with a larger distance, and we observe that these three methods behave similarly; in Test 3, we choose two circles with a smaller distance, then the solution of the level set method separates, but the solutions of the phase field method and the energy minimization method merge unless $\epsilon$ is very small, in which case the computational cost is large (some spatial grids should be placed iniside the diffuse interface); in Test 4, we choose another wedge-like initial condition, then the solution of the level set method merges, but the solutions of the phase field method and the energy minimization method separate unless $\epsilon$ is very small.\\%We also give the energy plots for all these cases.\\

\begin{figure}[H]
\centering
\includegraphics[height=2.4in,width=3.4in]{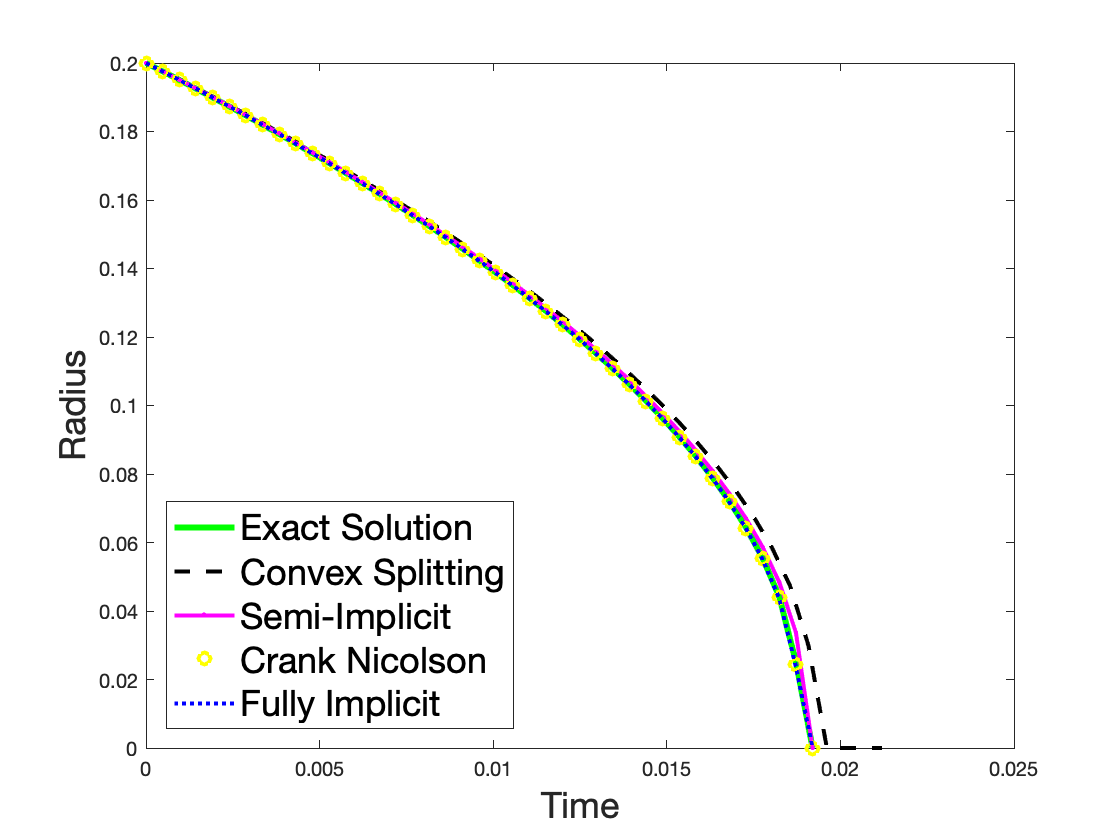}
\caption{{\color{black}Test 1.} Evolutions of the solutions of different numerical schemes}
\label{fig1}
\end{figure}

{\bf Test 1.} \label{test1} In this test, we test a benchmark problem with a circle with $radius = 0.2$ for the initial condition based on different numerical schemes. Figure \ref{fig1} shows the change of radius with respect to time for different numerical schemes. As expected, we find that when the time step size is sufficiently small, the evolutions of the solutions are very similar when different numerical schemes are used. Therefore, we just need to use the fully implicit scheme with a small time step size to compute the evolutions of the solutions.

\begin{figure}[H]
\centering
\includegraphics[height=2.0in,width=4.6in]{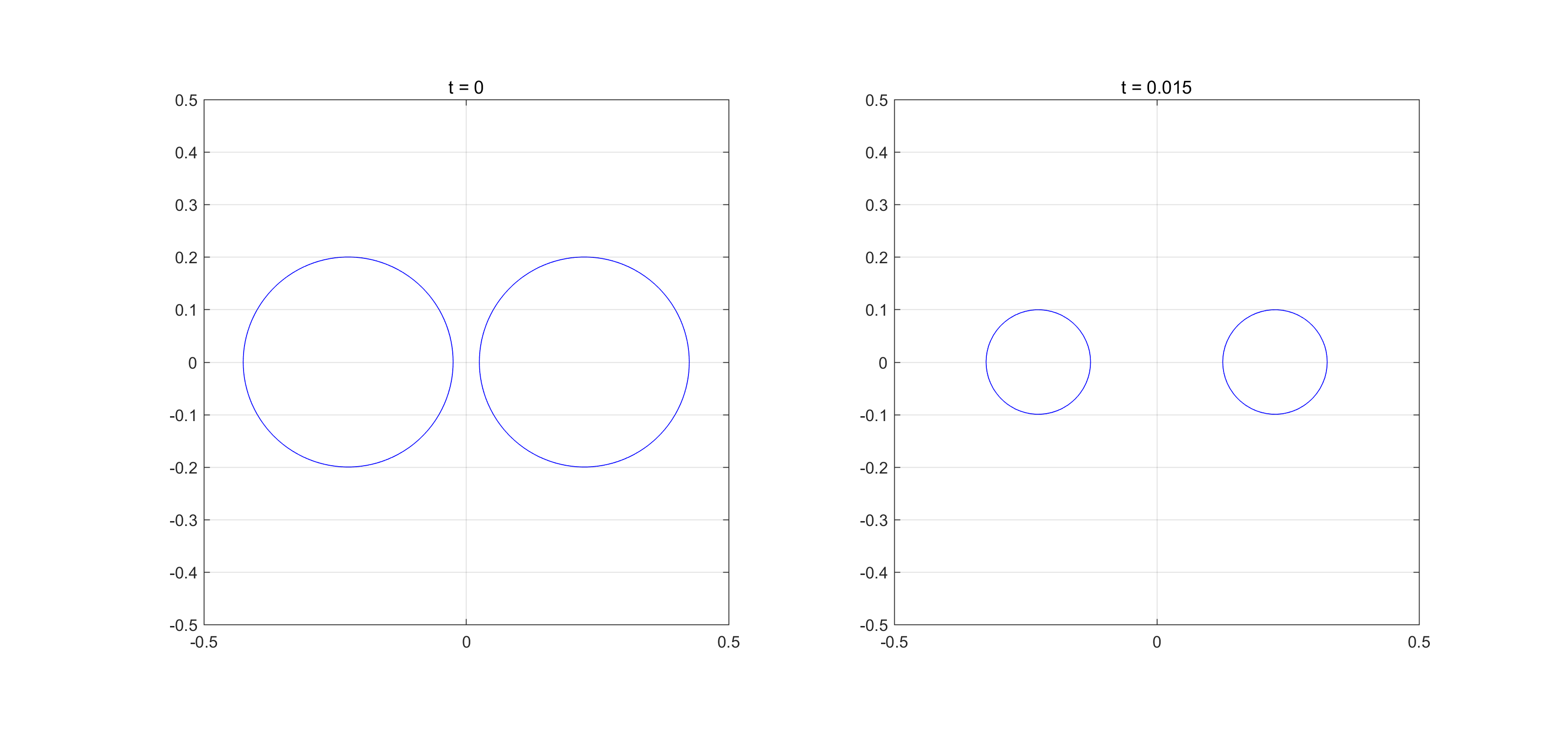}
\caption{{\color{black} Test 2. }The level set method is used. The left graph is the initial condition, and the right graph is the evolution of the solution at $t=0.015$.}
\label{fig2}
\end{figure}

\begin{figure}[H]
\centering
\includegraphics[height=1.4in,width=2.3in]{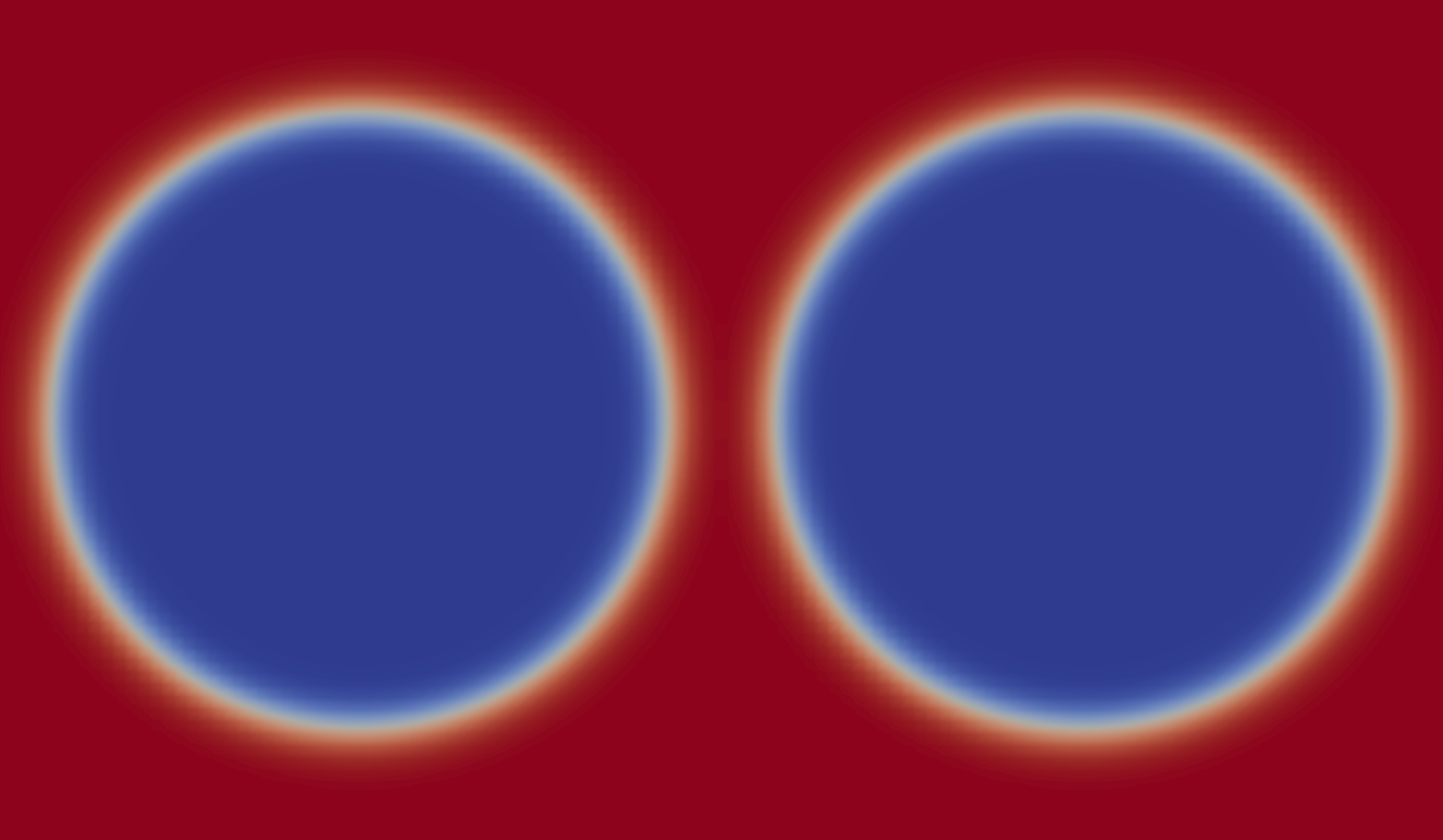}
\includegraphics[height=1.4in,width=2.3in]{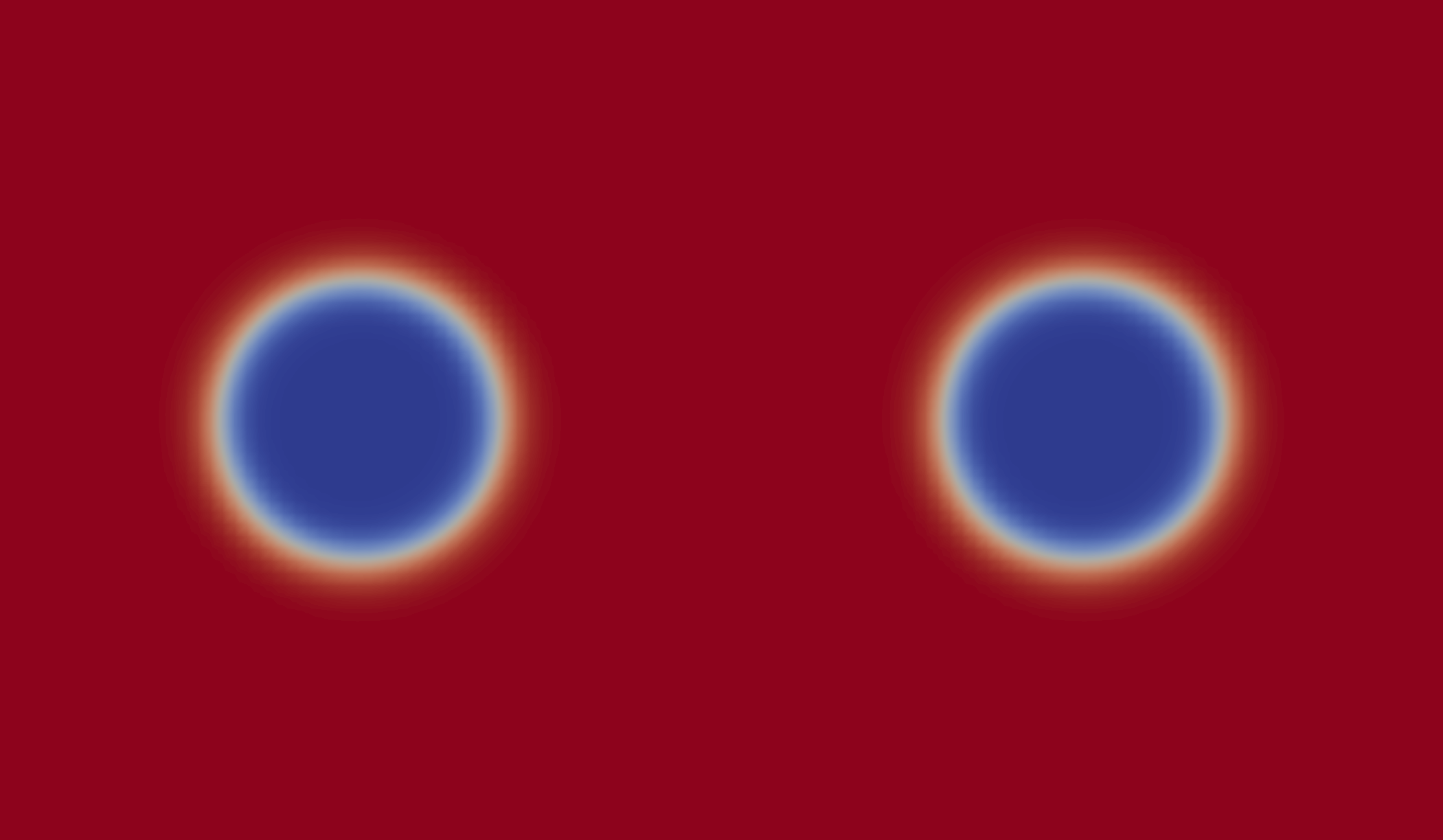}
\caption{{\color{black} Test 2. }The phase field method is used. The left graph is the initial condition, and the right graph is the evolution of the solution at $t=0.015$ {\color{black} with $\epsilon = 0.01$}.}
\label{fig3}
\end{figure}

\begin{figure}[H]
\centering
\includegraphics[height=1.4in,width=2.3in]{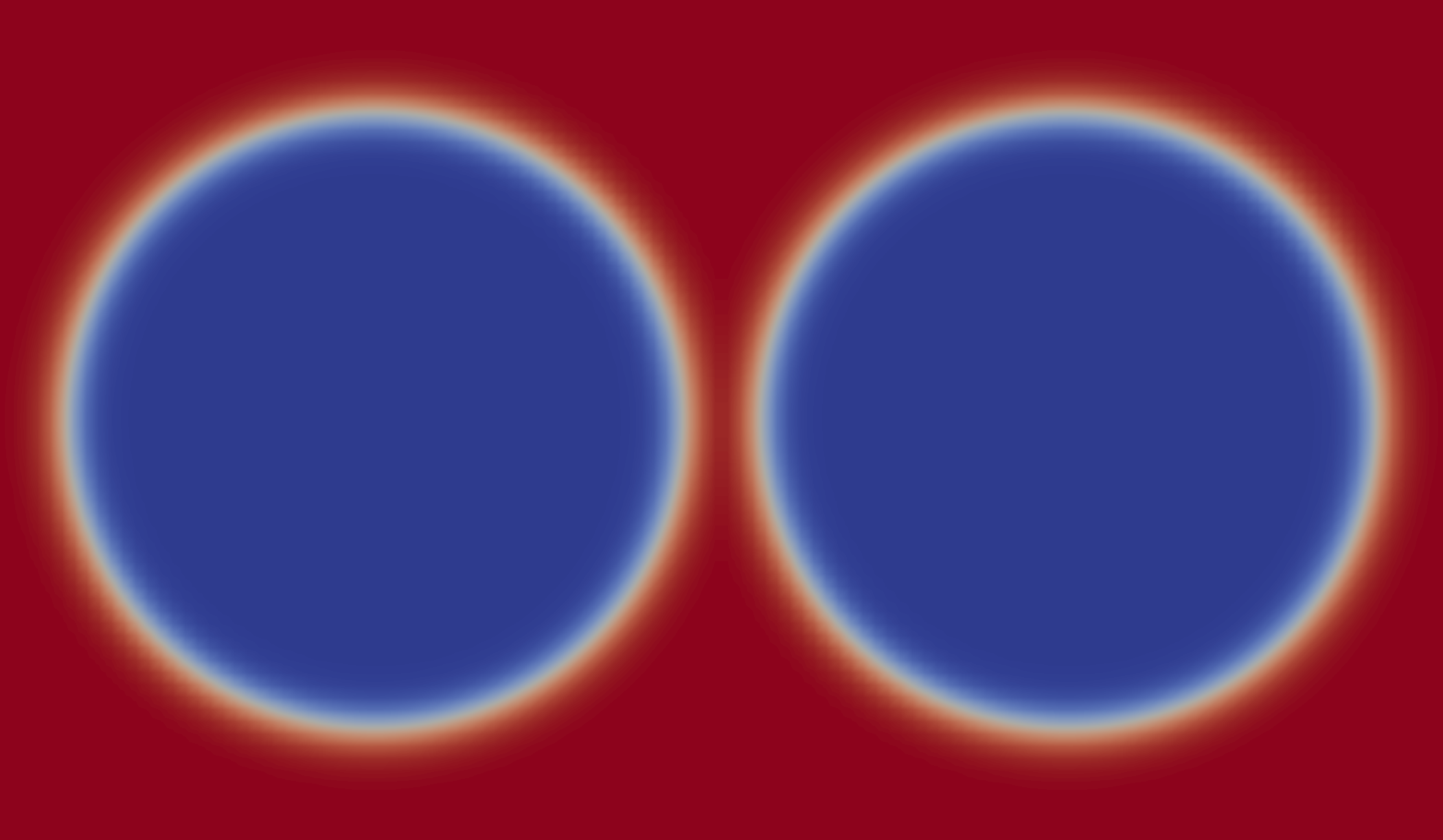}
\includegraphics[height=1.4in,width=2.3in]{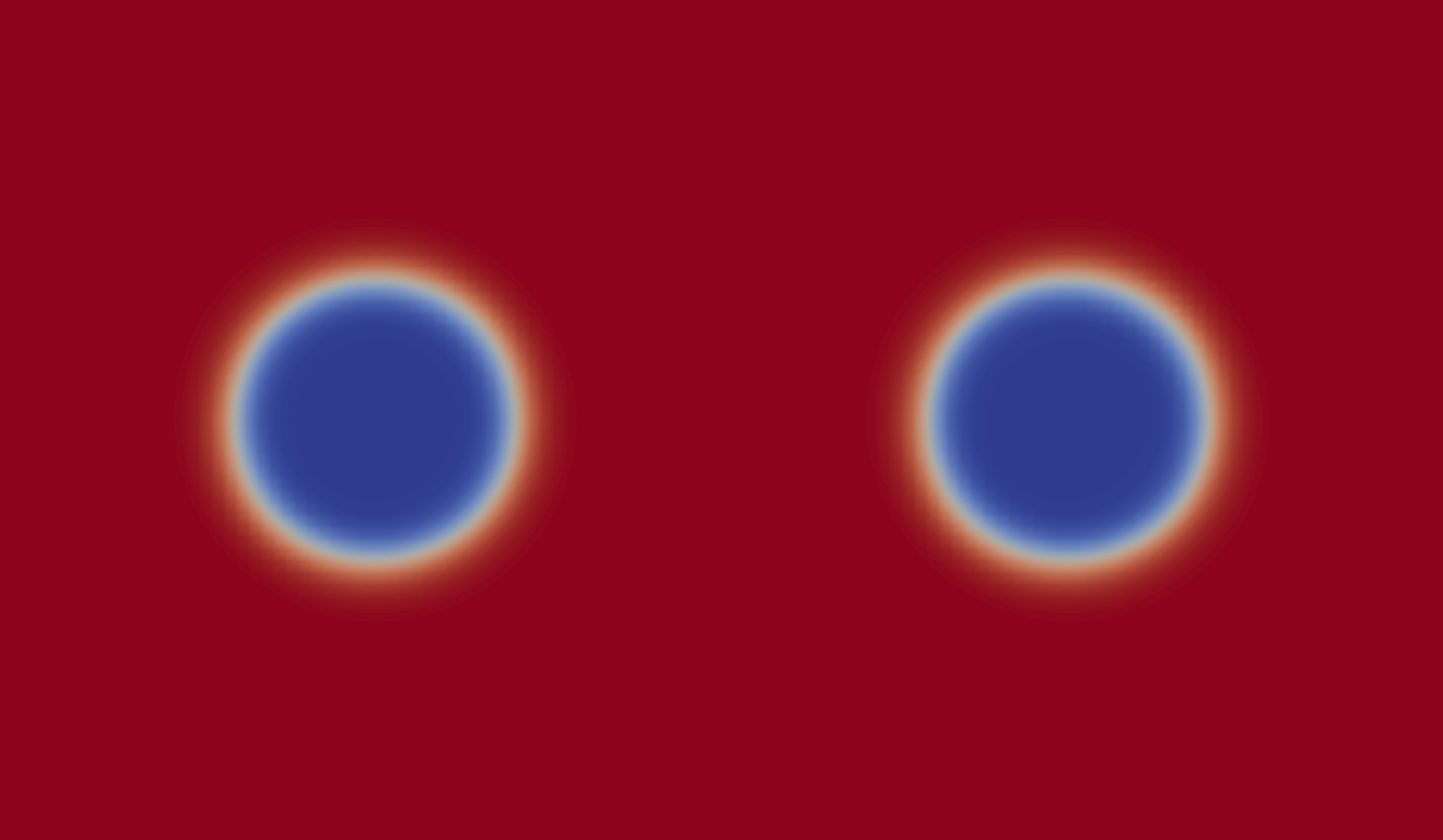}
\caption{{\color{black} Test 2. }The energy minimization method is used. The left graph is the initial condition, and the right graph is the evolution of the solution at $t=0.015$ {\color{black} with $\epsilon = 0.01$}.}
\label{fig4}
\end{figure}

\begin{figure}[H]
\centering
\includegraphics[height=2.2in,width=2.3in]{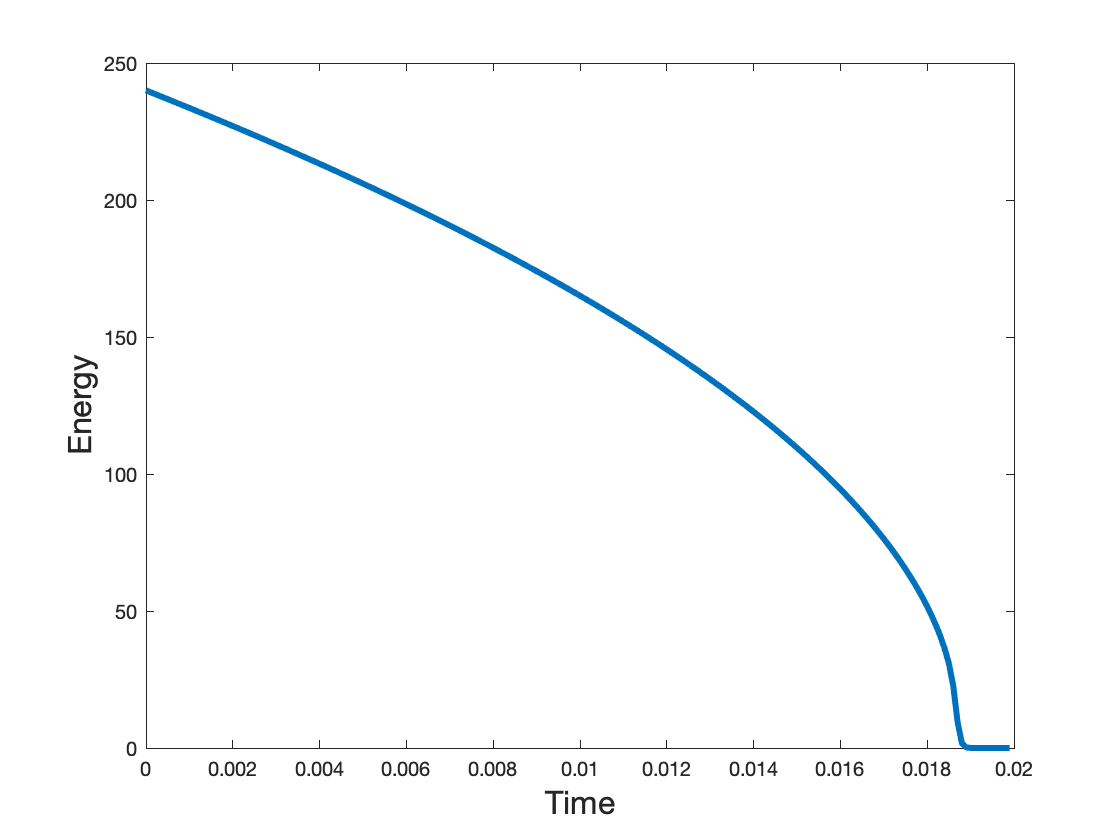}
\includegraphics[height=2.2in,width=2.3in]{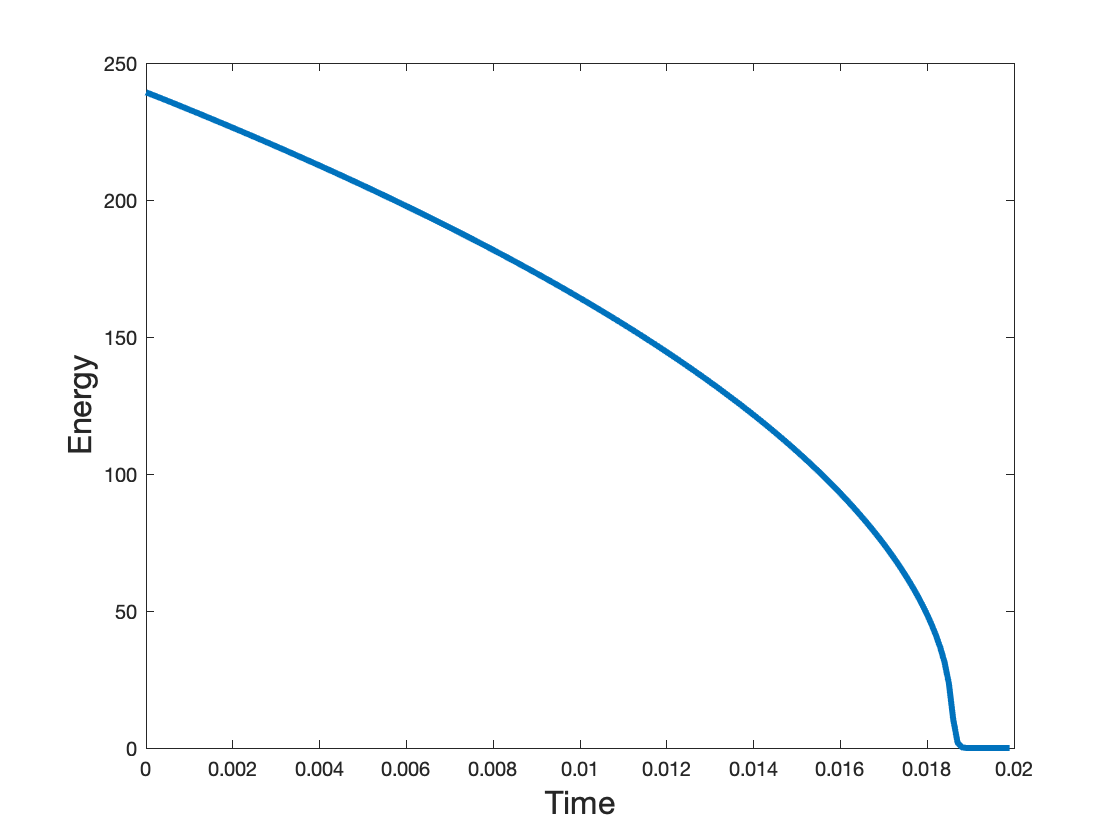}
\caption{Test 2. The left hand side is the energy of the FIS, and the right hand side is the energy of the energy minimization method with $\epsilon = 0.01$.}
\label{energy1}
\end{figure}

{\bf Test 2.} \label{test2} In this test, we compare the evolutions of the solutions when the initial condition is two circles with a large distance based on the level set method and the phase field method. In Figure \ref{fig2}, the level set method is used to compute the evolution of the solution, and we find that these two circles separate. The distance between these two circles is $d=0.05$, the spatial size is $h=0.01$, and the time step size is $k = 2.5\times10^{-5}$. In Figure \ref{fig3}, the fully implicit scheme of the phase field method is used to compute the evolution of the solution, and we find that these two circles separate. The distance between these two circles is $d=0.05$, the spatial size is $h=0.008$, the interaction length is $\epsilon=0.01$, and the time step size is $k = 1\times10^{-4}$. In Figure \ref{fig4}, {\color{black}using the same data for Figure \ref{fig3},} the energy minimization method is used to compute the evolution of the solution, and we find that these two circles separate. From this test, we observe that the evolutions of the solutions based on different methods behave similarly. The Figure \ref{energy1} indicates the energy change over time for both the phase field method and the energy minimization method. This test takes the initial condition used in many other papers, and this test demonstrates that the topological changes need not be considered when the distance between the two circles is large.
%; one may observe that both methods' energy decreases gradually until time 0.018, and then there's a sharp decrease until it reaches 0.  

%:

\begin{figure}[H]
\centering
\includegraphics[height=2.0in,width=4.6in]{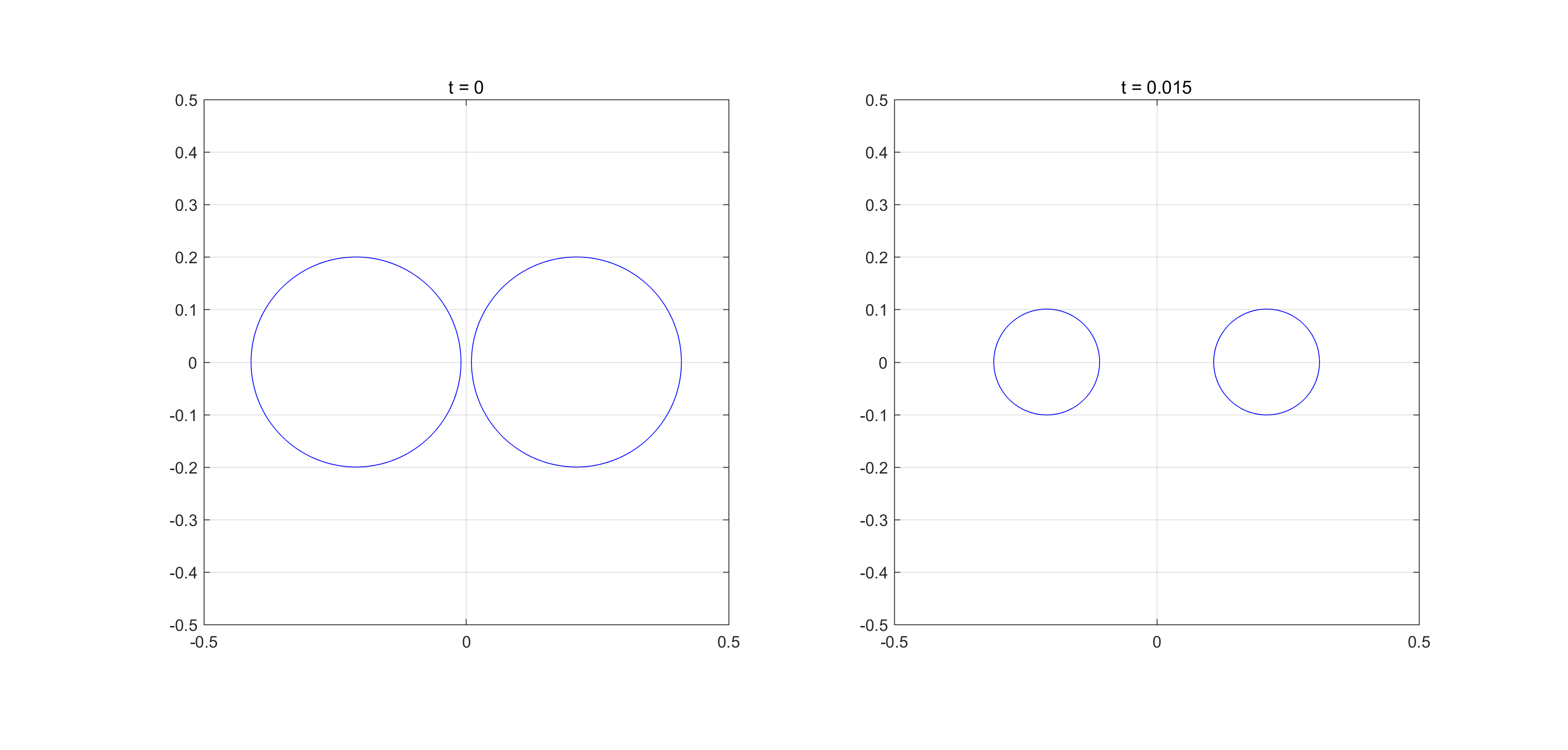}
\caption{{\color{black} Test 3. }The level set method is used. The left graph is the initial condition, and the right graph is the evolution of the solution at $t=0.015$.}
\label{fig5}
\end{figure}

\begin{figure}[H]
\centering
\includegraphics[height=1.4in,width=2.3in]{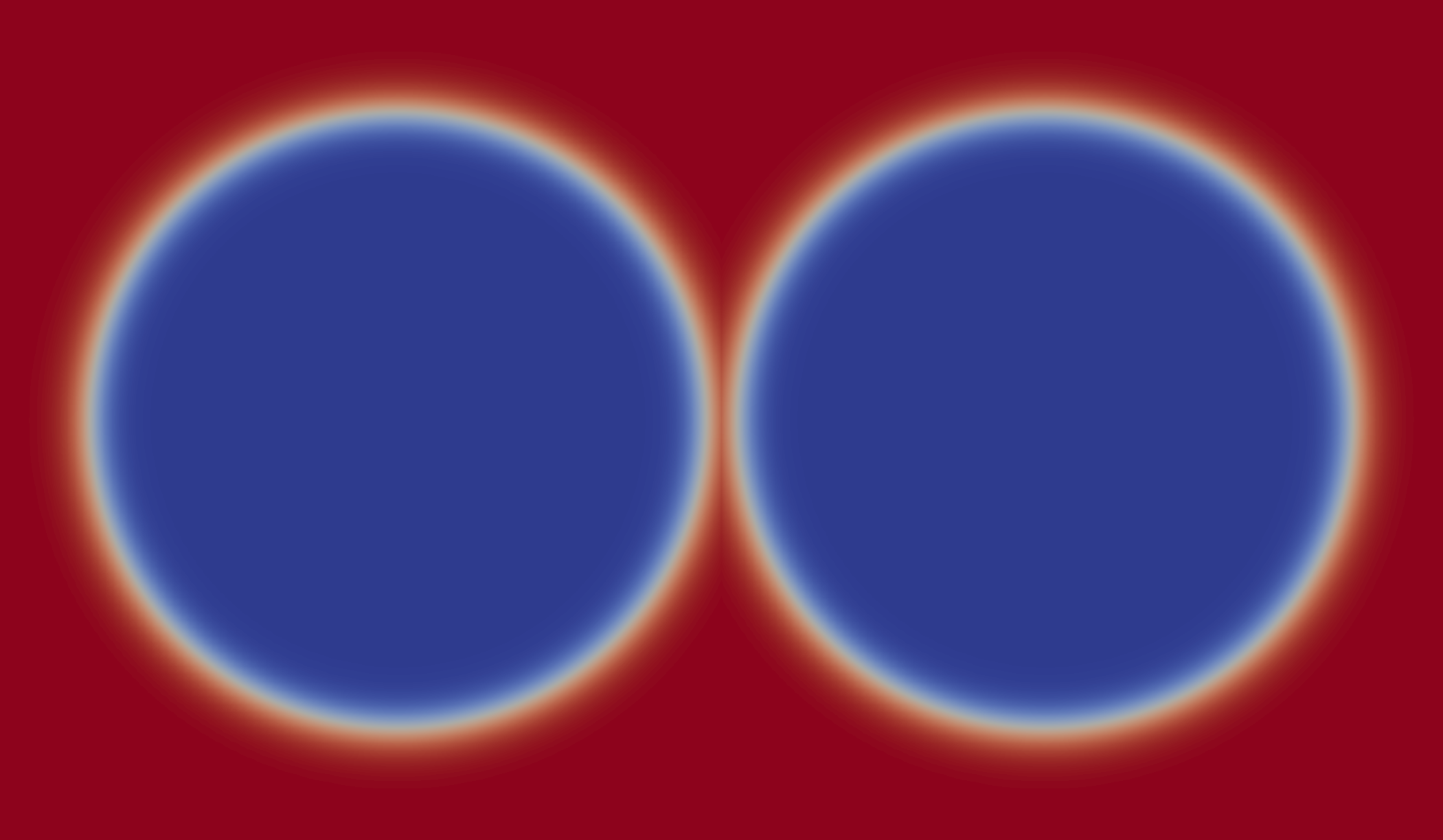}
\includegraphics[height=1.4in,width=2.3in]{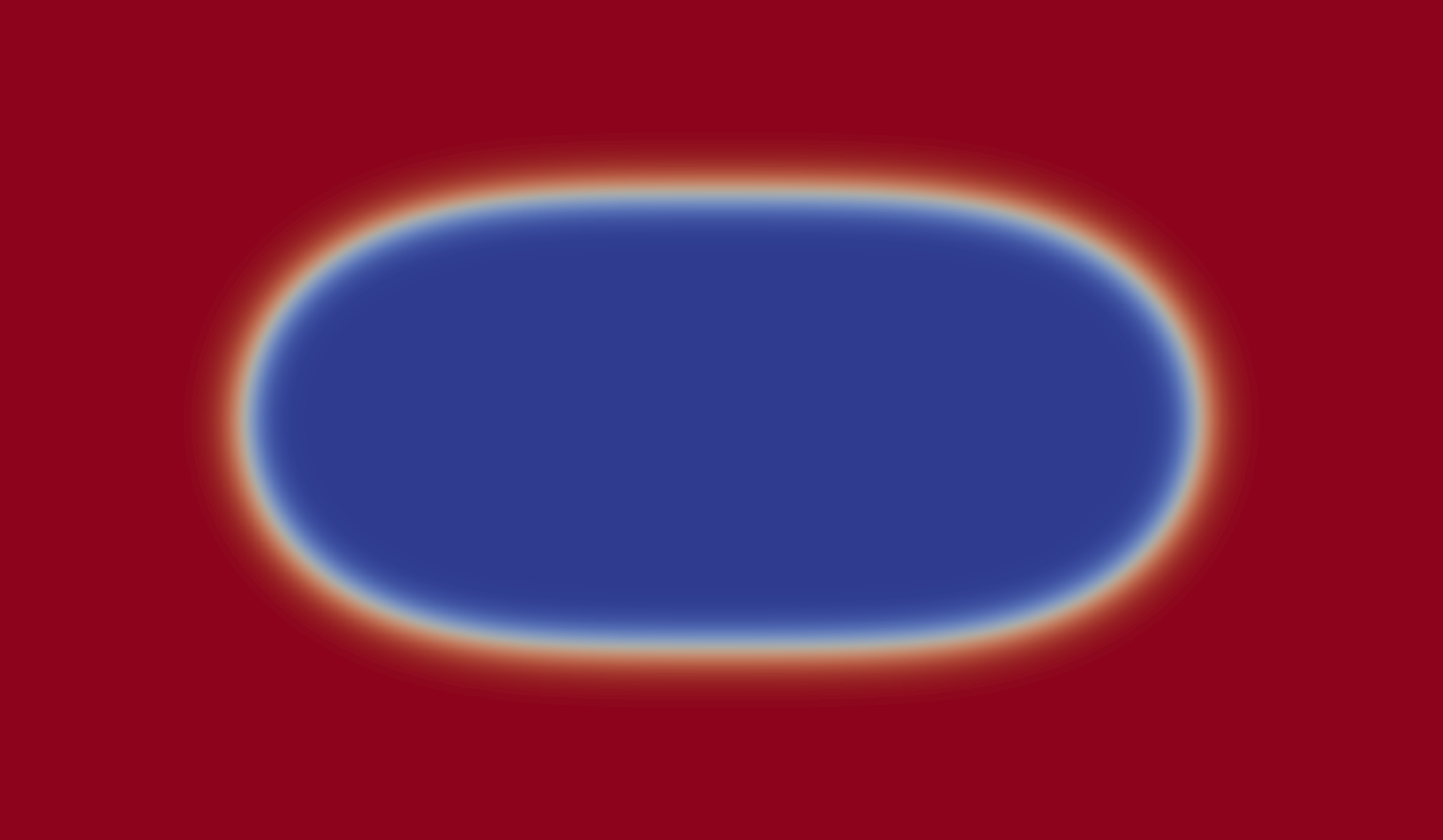}
\caption{{\color{black} Test 3. }The phase field method is used. The left graph is the initial condition, and the right graph is the evolution of the solution at $t=0.015$ {\color{black} with $\epsilon = 0.01$}.}
\label{fig6}
\end{figure}

\begin{figure}[H]
\centering
\includegraphics[height=1.4in,width=2.3in]{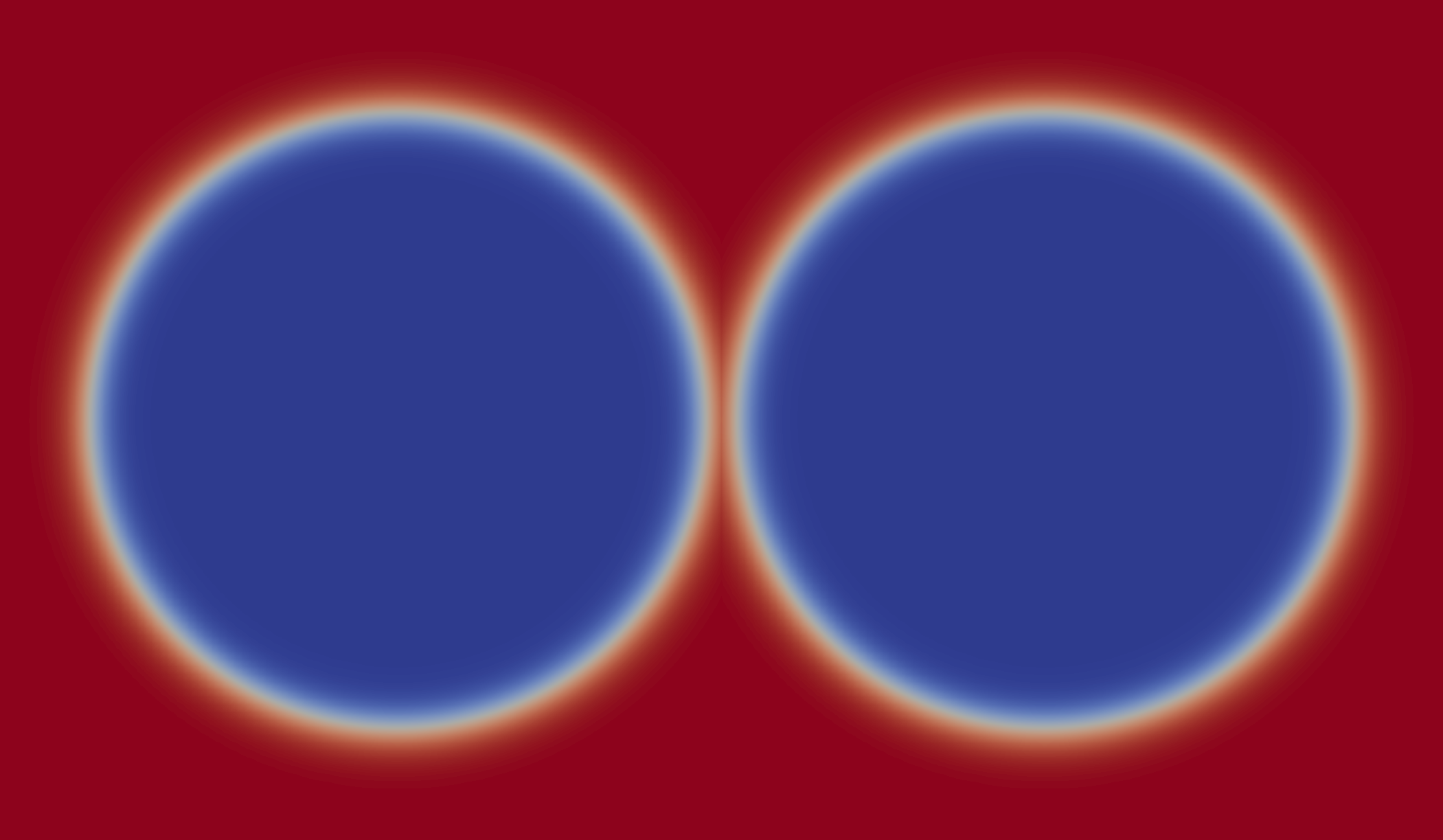}
\includegraphics[height=1.4in,width=2.3in]{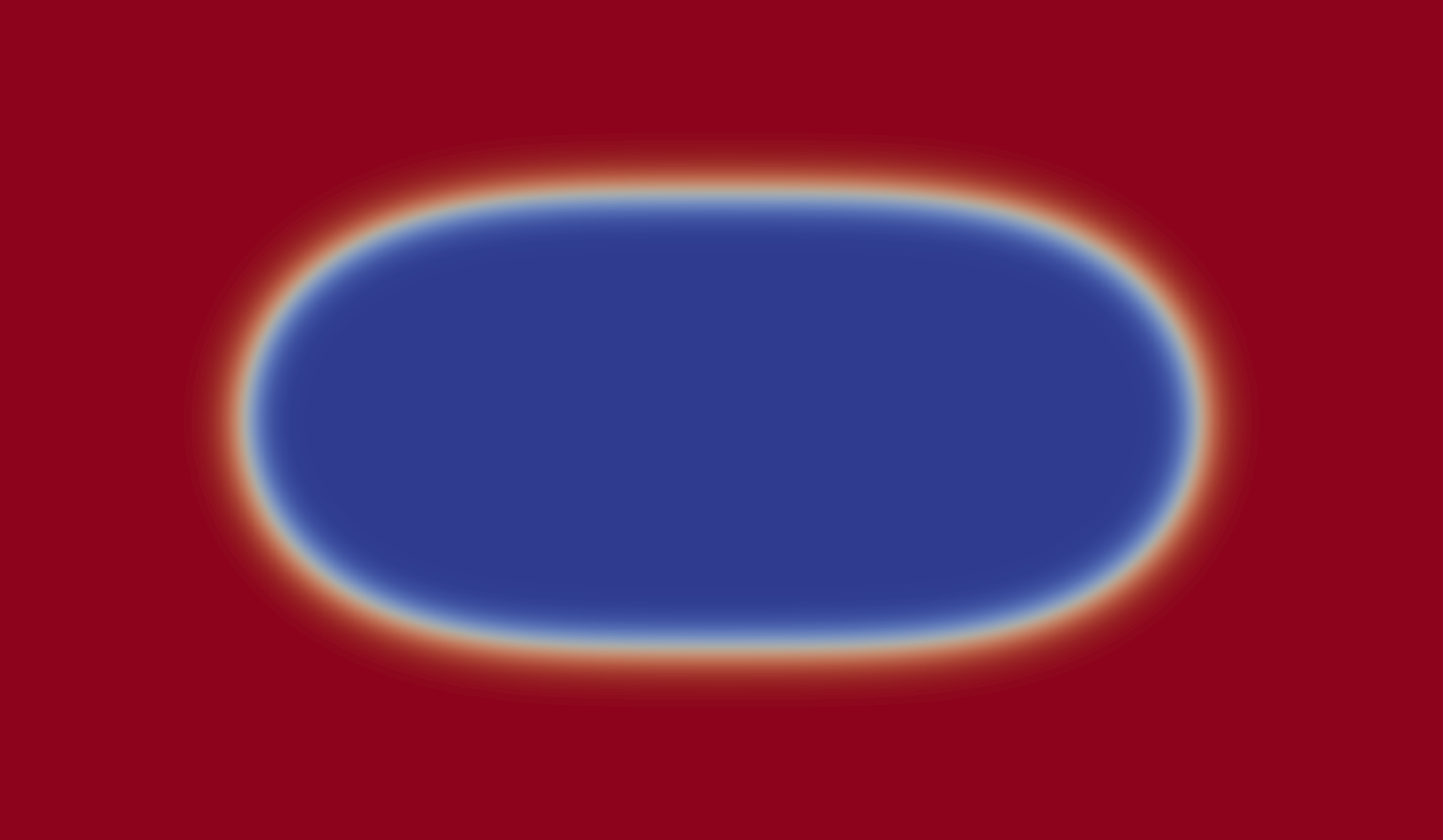}
\caption{{\color{black} Test 3. }The energy minimization method is used. The left graph is the initial condition, and the right graph is the evolution of the solution at $t=0.015$ {\color{black} with $\epsilon = 0.01$}.}
\label{fig7}
\end{figure}

\begin{figure}[H]
\centering
\includegraphics[height=2.2in,width=2.3in]{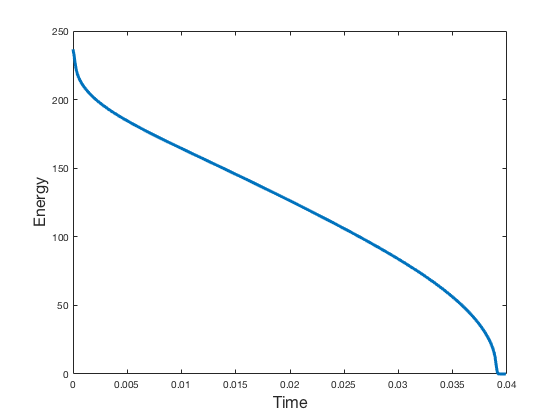}
\includegraphics[height=2.2in,width=2.3in]{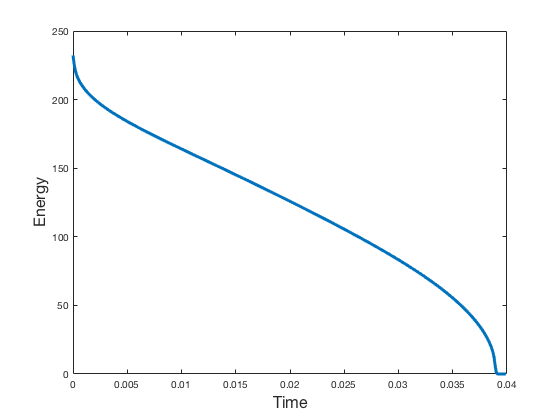}
\caption{Test 3. The left hand side is the energy of the FIS, and the right hand side is the energy of the energy minimization method with $\epsilon = 0.01$.}
\label{energy2}
\end{figure}

{\bf Test 3.} \label{test3} In this test, we compare the evolutions of the solutions when the initial condition is two circles with a small distance based on the level set method, the phase field method, and the energy minimization method. In Figure \ref{fig5}, the level set method is used to compute the evolution of the solution, and we find that these two circles separate. The distance between these two circles is $d=0.02$, the spatial size is $h=0.005$, and the time step size is $k = 2.5\times10^{-6}$. In Figure \ref{fig6}, the fully implicit scheme of the phase field method is used to compute the evolution of the solution, and we find that these two circles merge, which is inconsistent with the result from the level set method. The distance between these two circles is $d=0.02$, the spatial size is $h=0.005$, the interaction length is $\epsilon=0.01$, and the time step size is $k = 1\times10^{-4}$. In Figure \ref{fig7}, {\color{black}using the same data for Figure \ref{fig6},} the energy minimization method is used to compute the evolution of the solution, and we find that these two circles merge, which also contradicts the result of the level set method. Moreover, we tried much smaller time step sizes for the phase field method and the energy minimization method, but these two circles still merge. The Figure \ref{energy2} indicates the energy change over time for both the phase field method and the energy minimization method. We observe that the energy has a sharp decrease at the beginning, then it decreases slower until around time 0.035. Then, it decreases sharply again until it reaches 0. 

\begin{figure}[H]
\centering
\includegraphics[height=1.4in,width=2.3in]{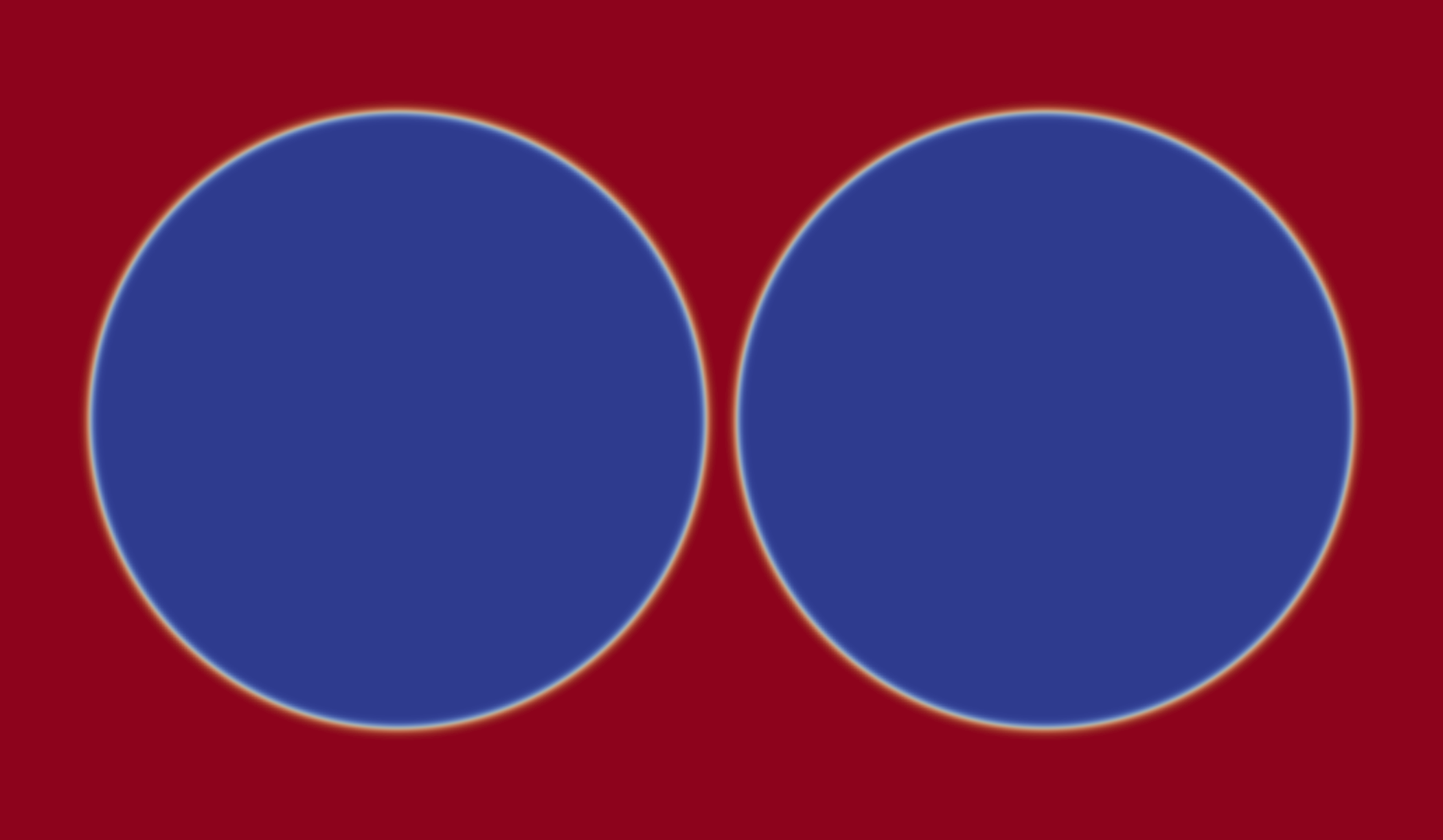}
\includegraphics[height=1.4in,width=2.3in]{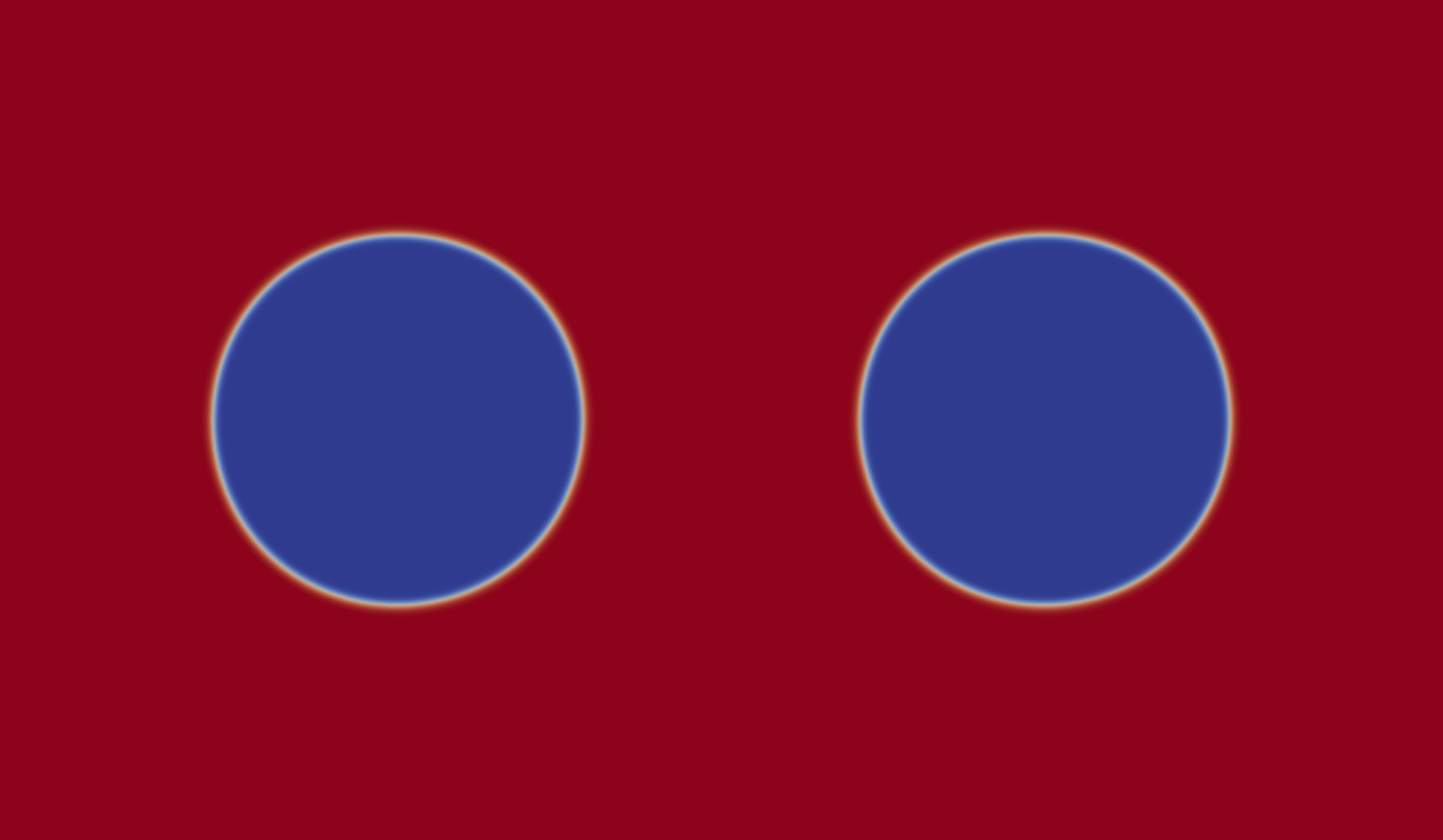}
\caption{{\color{black} Test 3. }The phase field method is used. The left graph is the initial condition, and the right graph is the evolution of the solution at $t=0.012$ {\color{black} with $\epsilon = 0.002$}.}
\label{fig8}
\end{figure}

\begin{figure}[H]
\centering
\includegraphics[height=1.4in,width=2.3in]{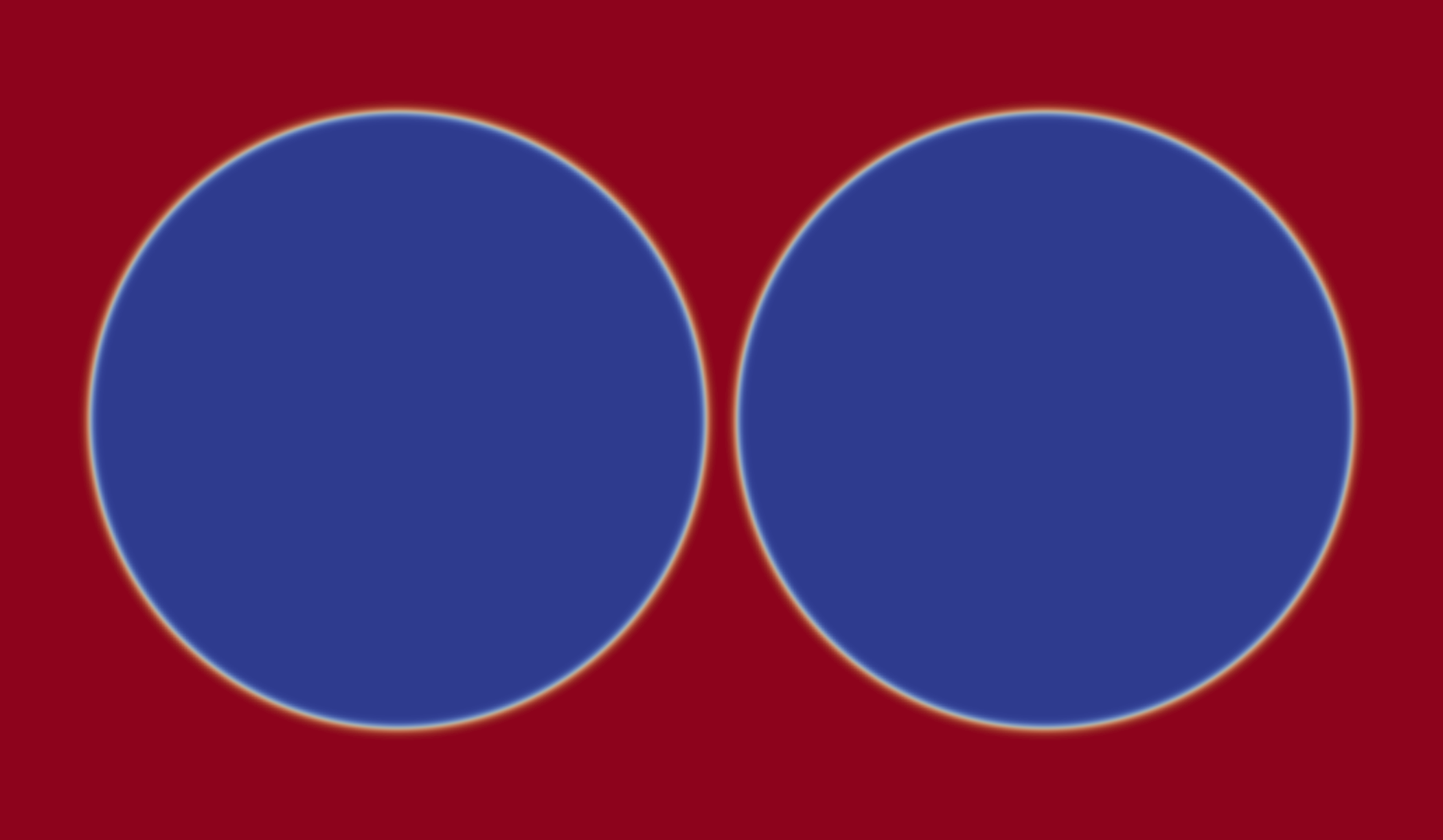}
\includegraphics[height=1.4in,width=2.3in]{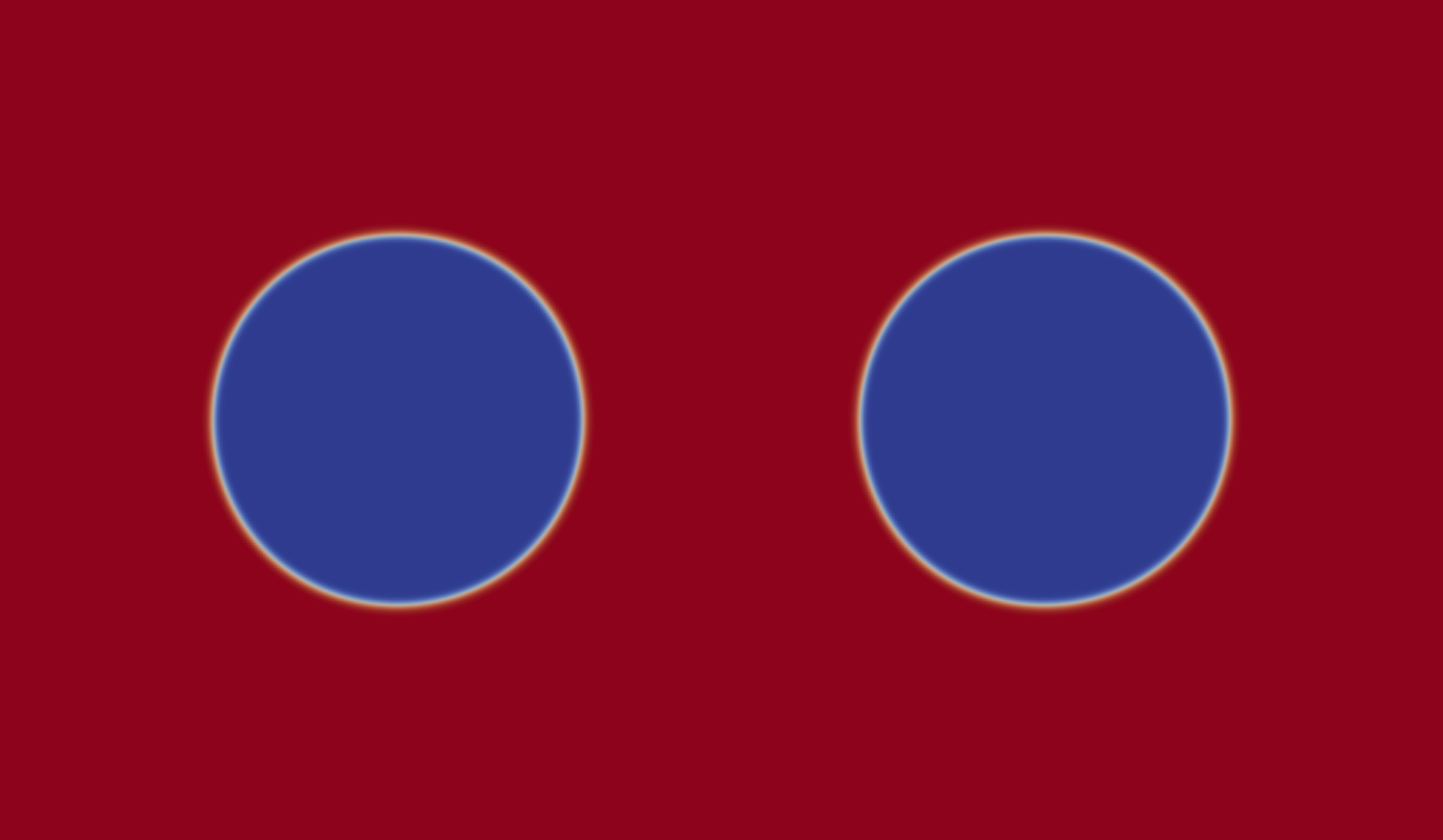}
\caption{{\color{black} Test 3. }The energy minimization method is used. The left graph is the initial condition, and the right graph is the evolution of the solution at $t=0.012$ {\color{black} with $\epsilon = 0.002$}.}
\label{fig9}
\end{figure}

\begin{figure}[H]
\centering
\includegraphics[height=2.2in,width=2.3in]{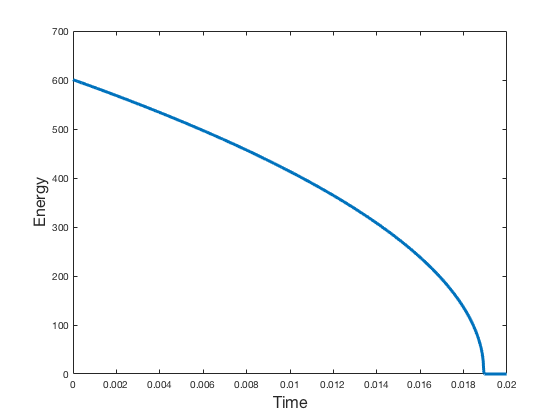}
\includegraphics[height=2.2in,width=2.3in]{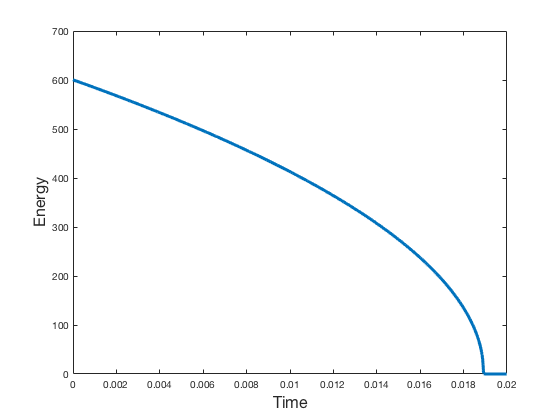}
\caption{Test 3. The left hand side is the energy of the FIS, and the right hand side is the energy of the energy minimization method with $\epsilon = 0.002$.}
\label{energy3}
\end{figure}

In Figure \ref{fig8}, the phase field method is used to compute the evolution of the solution, and we find that these two circles separate. The distance between these two circles is $d=0.02$, the spatial size is $h=0.0018$, the interaction length is $\epsilon=0.002$, and the time step size is $k = 4\times10^{-6}$. In Figure \ref{fig9}, {\color{black}using the same data for Figure \ref{fig8},} the energy minimization method is used to compute the evolution of the solution, and we find that these two circles separate. Using a small $\epsilon$, those two methods give very similar results as the level set method. However, with a small $\epsilon$, the required computational cost increases significantly. The Figure \ref{energy3} indicates the energy change over time for both the phase field method and the energy minimization method. We observe that the energy decreases gradually until time 0.018, and then decreases sharply until it reaches 0. This test validates that the interaction length plays a crucial role for some initial conditions. 

%In Figure \ref{fig10}, the energy minimization method is used to compute the evolutions, and we find that these two circles will separate. Here the distance between these two circle $d=0.02$, spatial size $h=\frac{1}{550}$, interaction length $\epsilon=0.002$, and time step size $k = 4\times10^{-5}$, which is $k = 100\times\epsilon^2$. Phase field method will not converge if $k = 100\times\epsilon^2$. 
%
%\begin{figure}[H]
%\centering
%\includegraphics[height=1.5in,width=2.5in]{pic/test3_pic10.png}
%\includegraphics[height=1.5in,width=2.5in]{pic/test3_pic11.png}
%\caption{The energy minimization method is used. The left graph is the initial condition, and the right graph is the evolution at $t=0.0024$.}
%\label{fig10}
%\end{figure}

\begin{figure}[H]
\centering
\includegraphics[height=2.2in,width=4.6in]{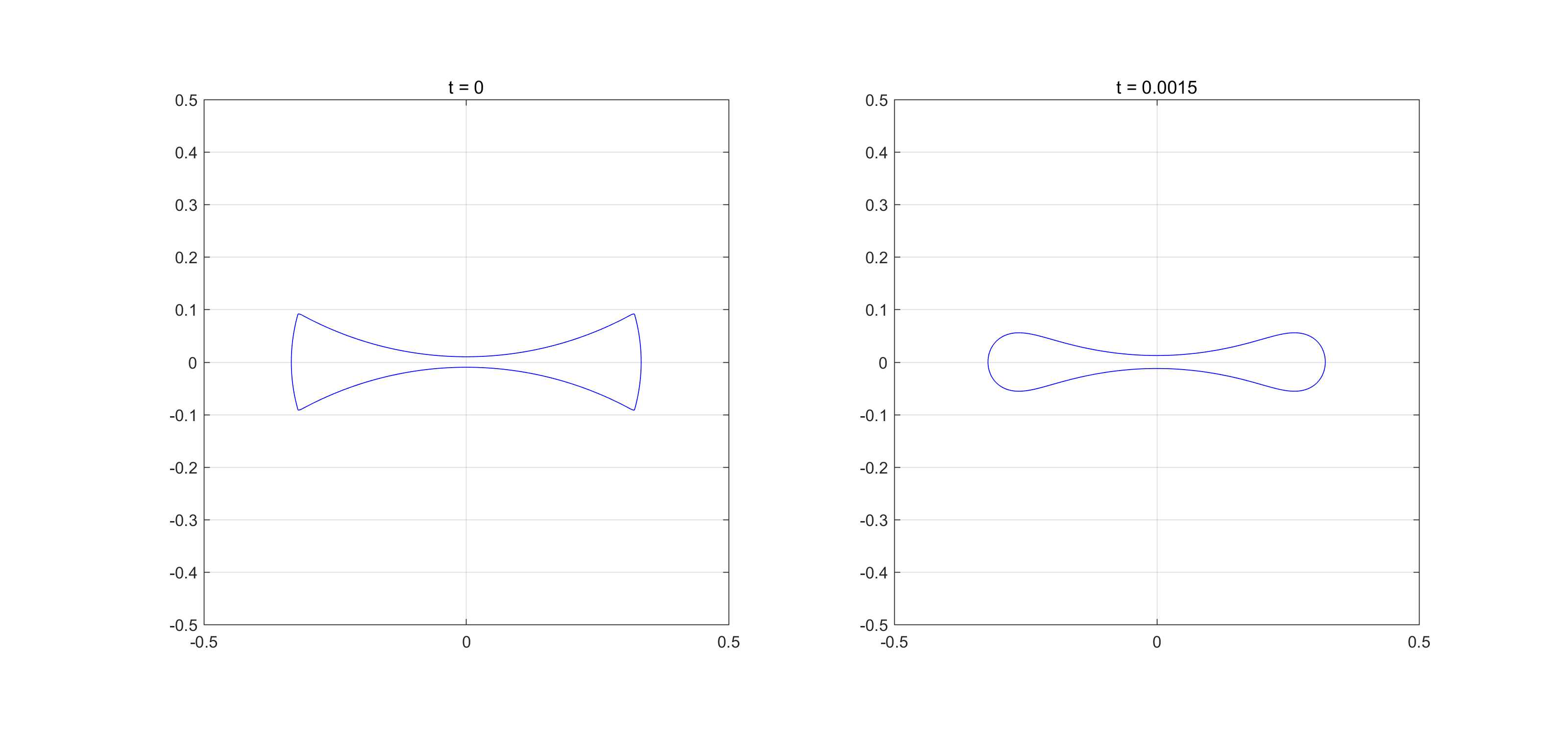}
\caption{{\color{black} Test 4. }The level set method is used. The left graph is the initial condition, and the right graph is the evolution of the solution at $t=0.015$.}
\label{fig11}
\end{figure}

\begin{figure}[H]
\centering
\includegraphics[height=1.4in,width=2.3in]{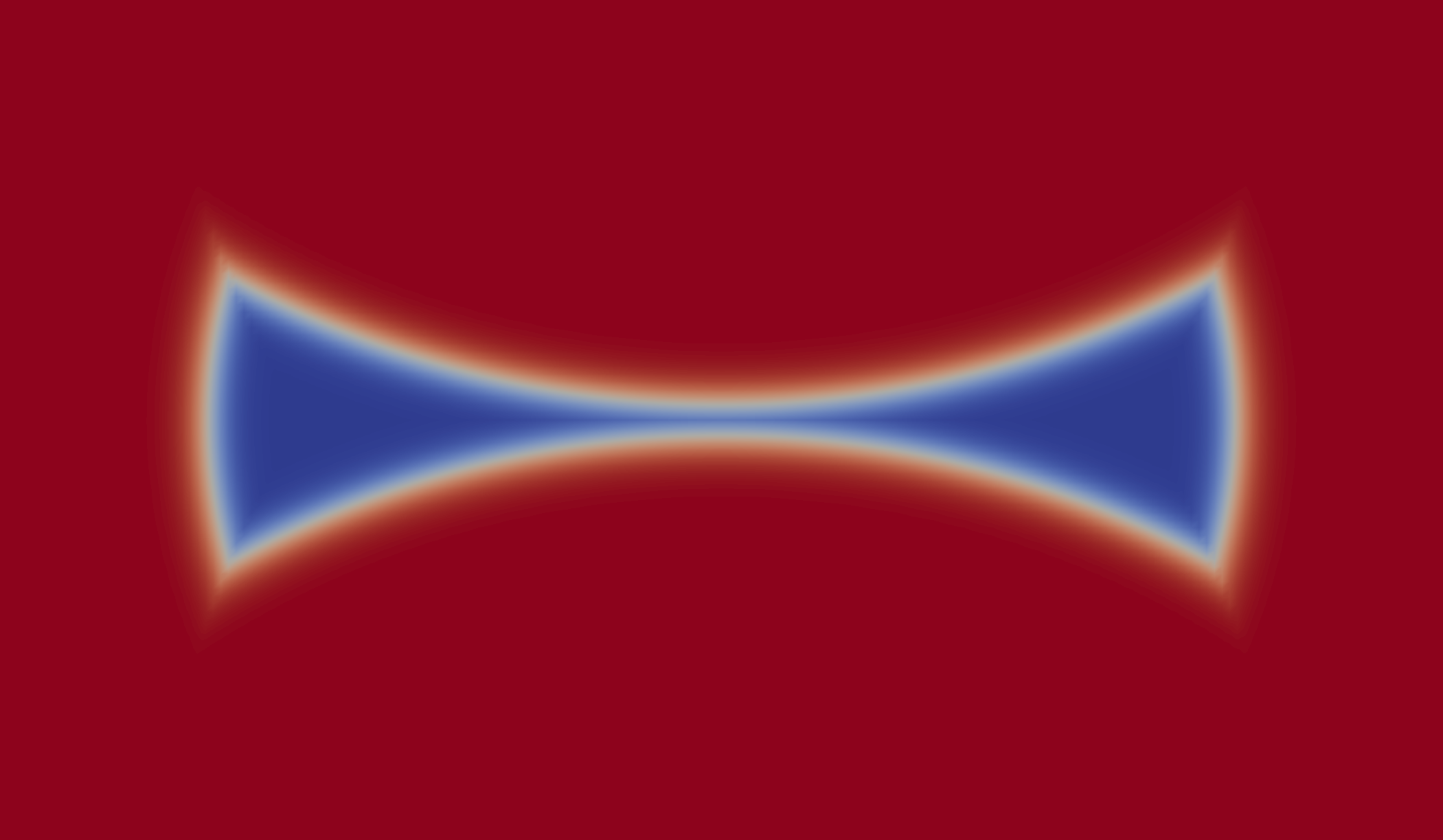}
\includegraphics[height=1.4in,width=2.3in]{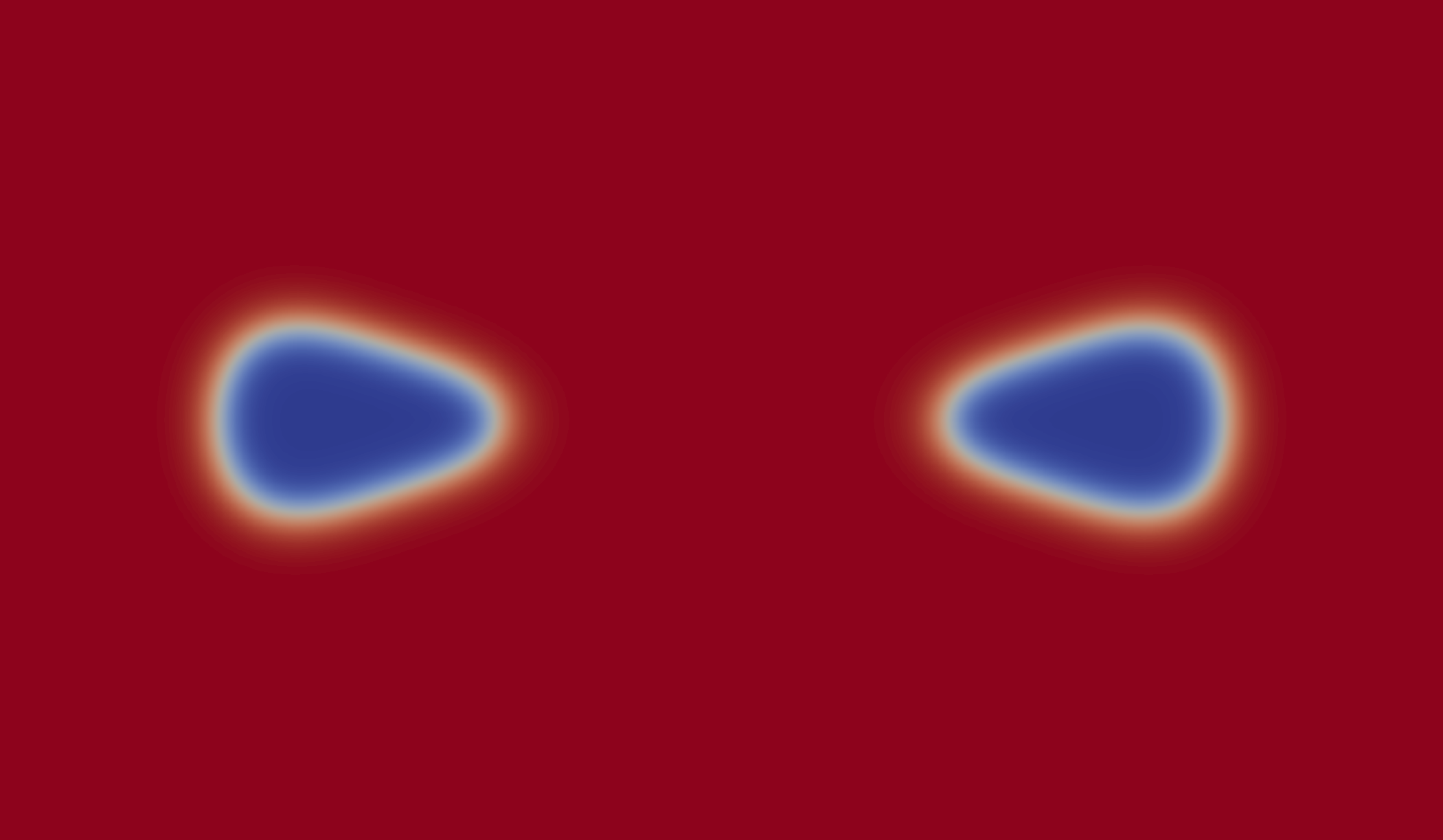}
\caption{{\color{black} Test 4. }The phase field method is used. The left graph is the initial condition, and the right graph is the evolution of the solution at $t=0.001$ {\color{black} with $\epsilon = 0.01$}.}
\label{fig12}
\end{figure}

\begin{figure}[H]
\centering
\includegraphics[height=1.4in,width=2.3in]{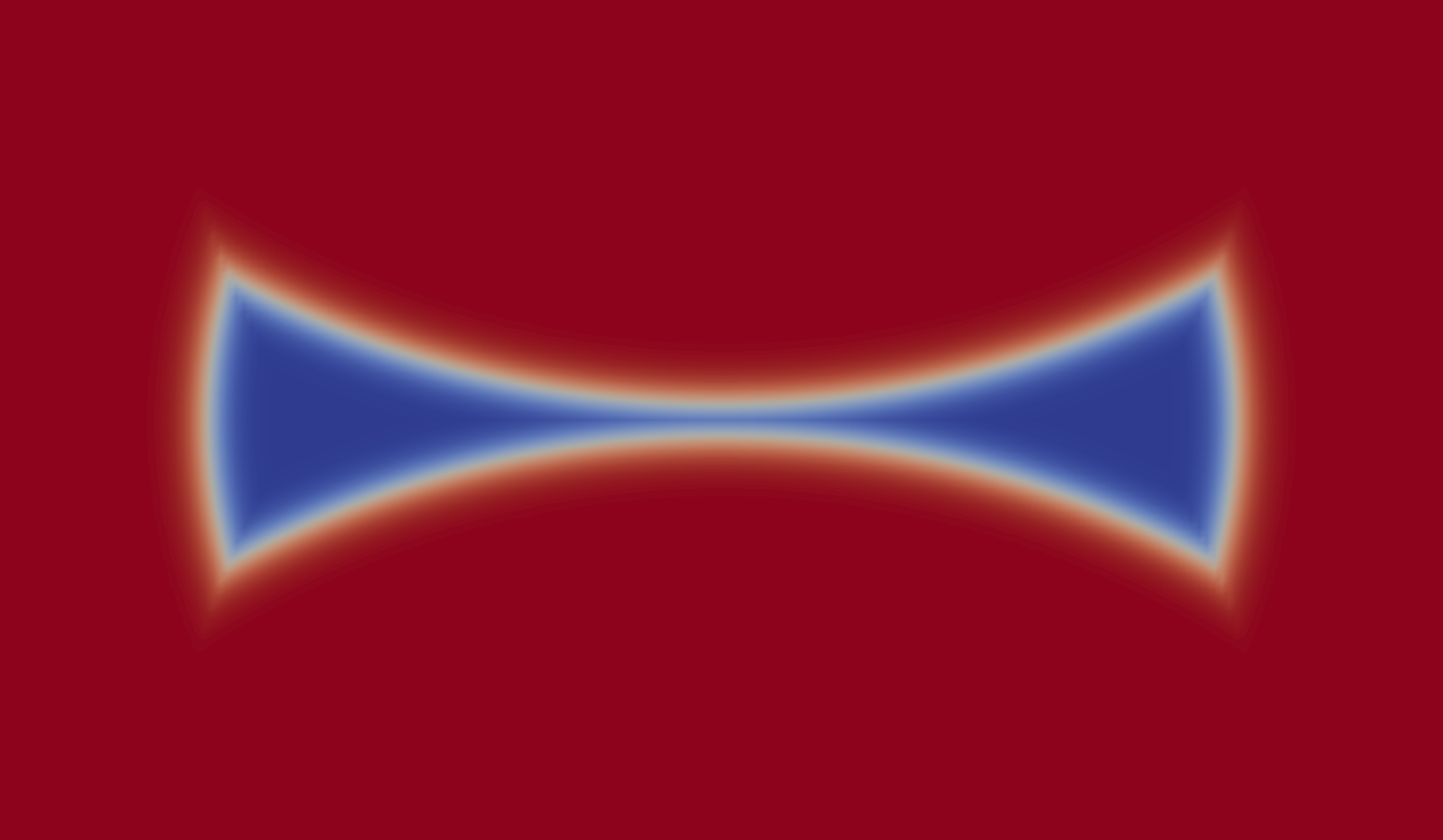}
\includegraphics[height=1.4in,width=2.3in]{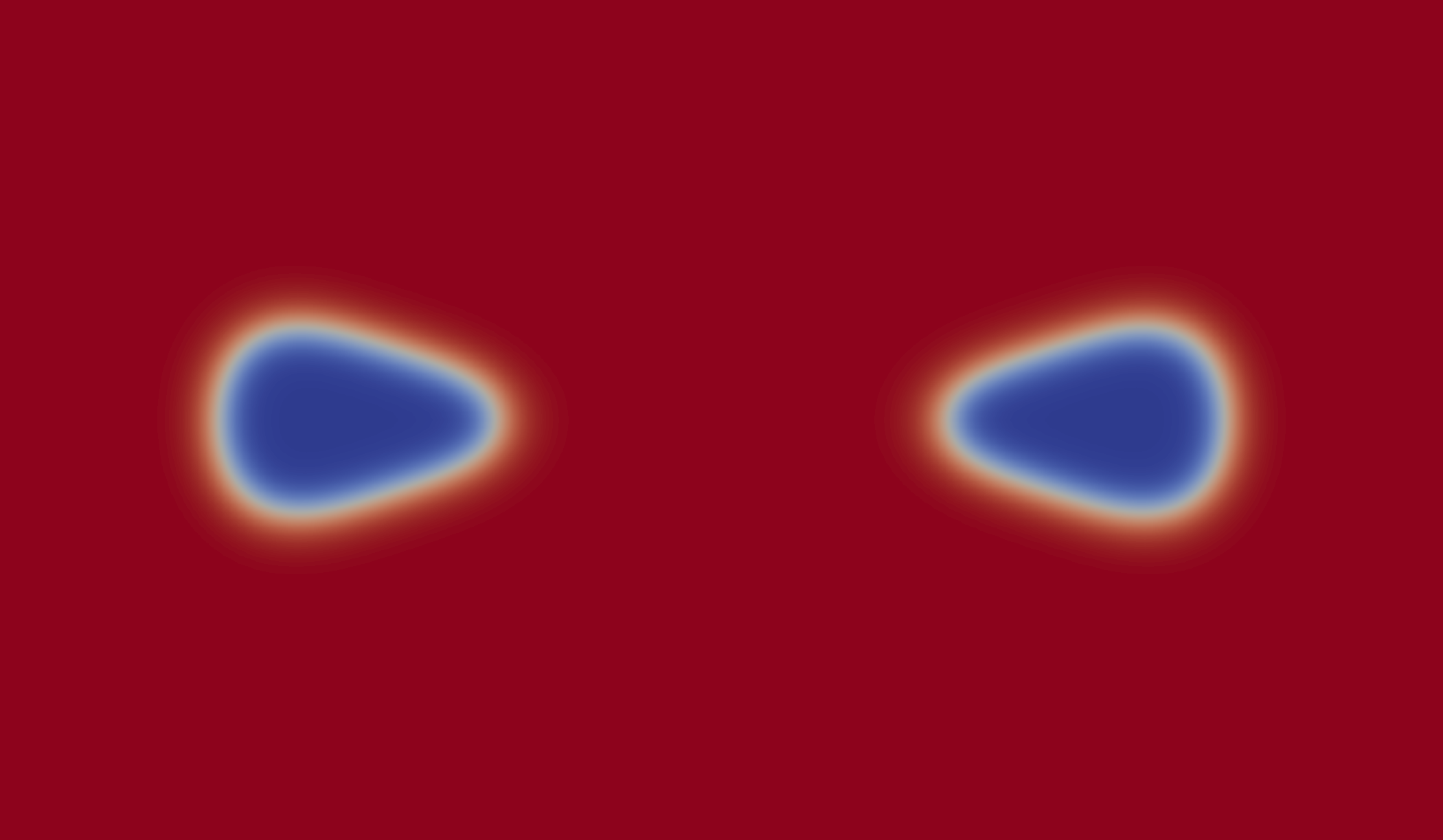}
\caption{{\color{black} Test 4. }The energy minimization method is used. The left graph is the initial condition, and the right graph is the evolution of the solution at $t=0.001$ {\color{black} with $\epsilon = 0.01$}.}
\label{fig13}
\end{figure}

\begin{figure}[H]
\centering
\includegraphics[height=2.2in,width=2.3in]{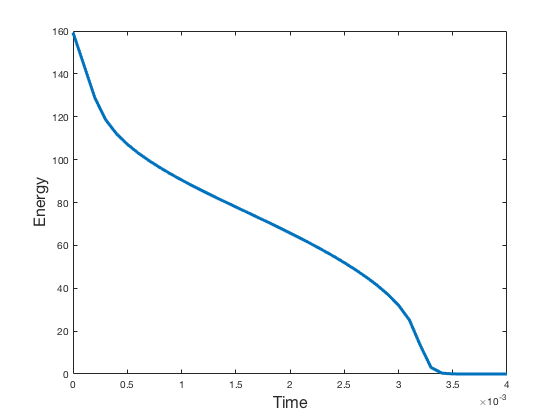}
\includegraphics[height=2.2in,width=2.3in]{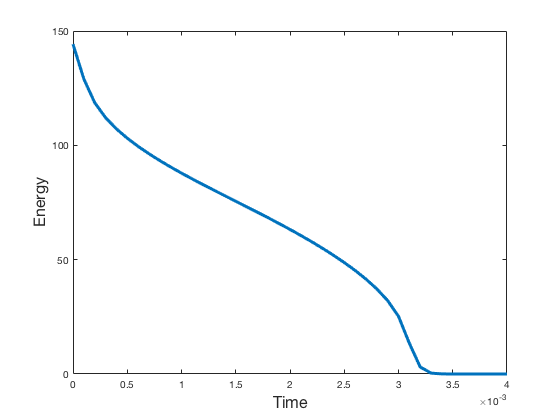}
\caption{{\color{black} Test 4. }The left hand side is the energy of the FIS, and the right hand side is the energy of the energy minimization method with $\epsilon = 0.01$.}
\label{energy5}
\end{figure}

{\bf Test 4.} %Check the wedge initial condition.
\label{test5} In this test, we compare the evolutions of the solutions when the initial condition is two wedges with small distance connections, based on the level set method, the phase field method, and the energy minimization method. Here the domain $\Omega: = [-0.5, 0.5]^2$. Define ${\bf m_1} = -{\bf m_3} := (0, 2)$, ${\bf m_2} := (0,0)$, $M := 0.01$, $r_1 = r_3 := 2/3 - 0.5\times M$, and $r_2 := 1/3$. Also, set $d_j({\bf X}) := |{\bf X} - {\bf m_j}| - r_j$ for ${\bf X} \in \Omega$, where ${\bf X} = (x, y)$. Then for ${\bf X} \in \Omega$, define
\begin{equation*}
d({\bf X}) := \max\{-d_1({\bf X}), d_2({\bf X}), -d_3({\bf X})\}, \hspace{0.6cm} u_0({\bf X}) := -\tanh(d({\bf X})/\sqrt{2}\times M).
\end{equation*}
Figure \ref{fig11} shows the evolution of the solution under the above initial condition using the level set method, with spatial size $h = 0.005$ and time step $k = 6.25\times10^{-6}$. We observe that the two wedges merge using the level set method. In Figure \ref{fig12}, the phase field method is used to compute the evolution of the solution, and we find that these two wedges separate, which is inconsistent with the result of the level set method. The spatial size is $h=0.005$, the interaction length is $\epsilon=0.01$, and the time step size is $k = 1\times10^{-4}$. In Figure \ref{fig13}, {\color{black}using the same data for Figure \ref{fig12},} the energy minimization method is used to compute the evolution of the solution, and we find that these two wedges separate, which is also inconsistent with the result of the level set method. Figure \ref{energy5} indicates the energy change over time for both the phase field method and the energy minimization method. We observe that the energy has a sharp drop at the beginning for both methods, then it decreases slowly until around time $3\times10^{-3}$. Then it decreases faster again until it reaches 0.   

\begin{figure}[H]
\centering
\includegraphics[height=1.4in,width=2.3in]{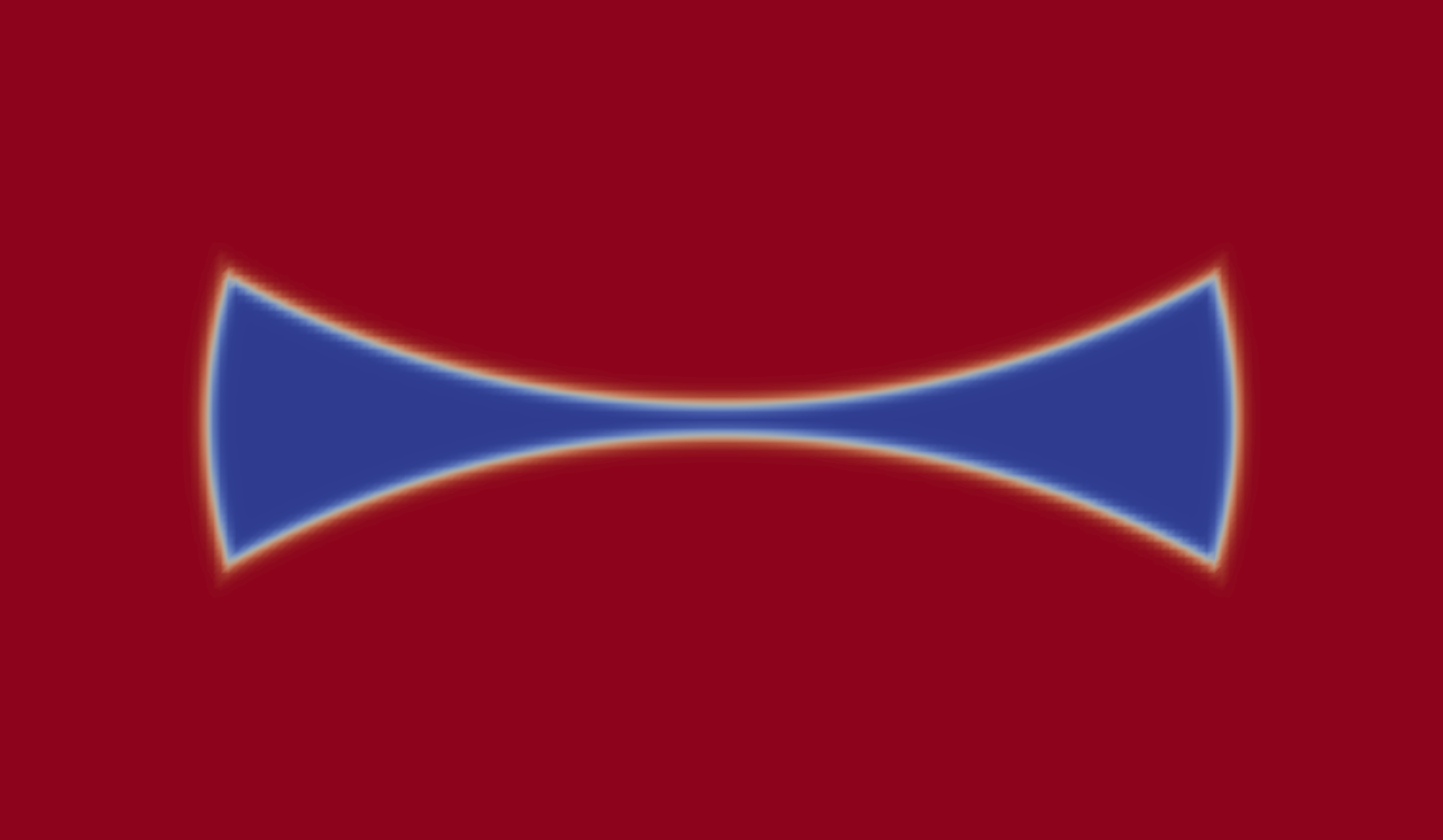}
\includegraphics[height=1.4in,width=2.3in]{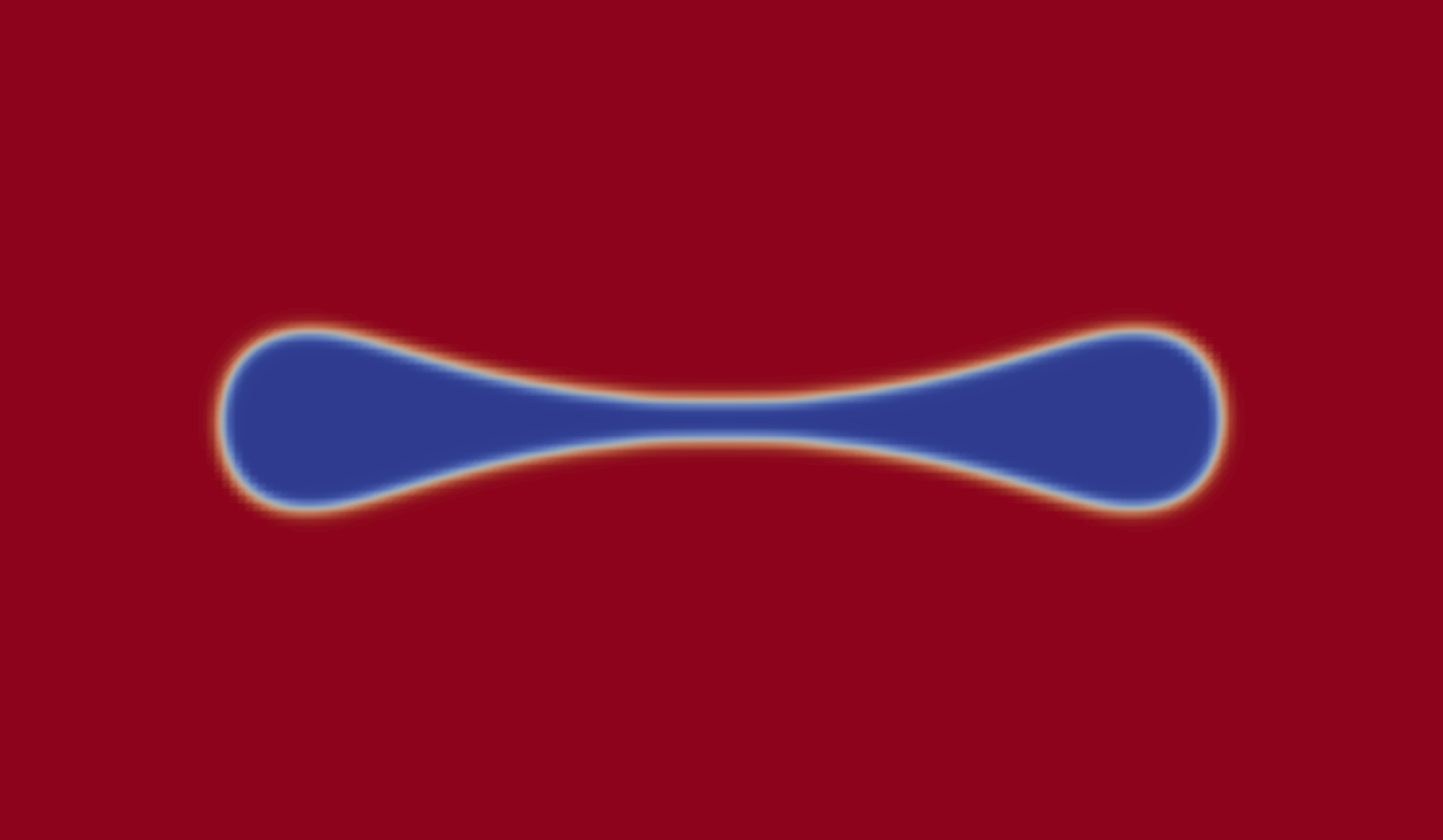}
\caption{{\color{black} Test 4. }The phase field method is used. The left graph is the initial condition, and the right graph is the evolution of the solution at around $t=0.00111$ {\color{black} with $\epsilon = 0.0033$}.}
\label{fig14}
\end{figure}

\begin{figure}[H]
\centering
\includegraphics[height=1.4in,width=2.3in]{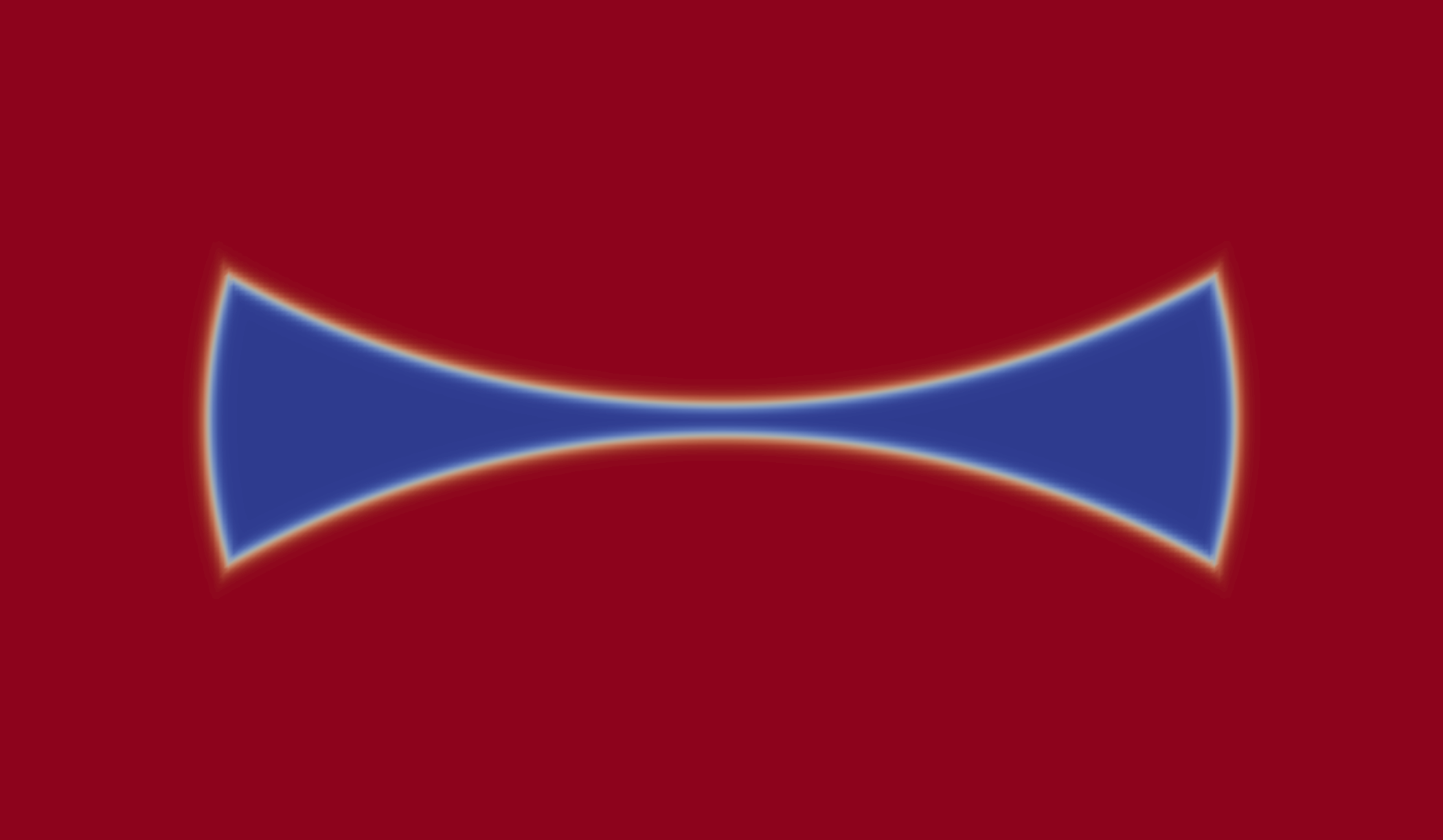}
\includegraphics[height=1.4in,width=2.3in]{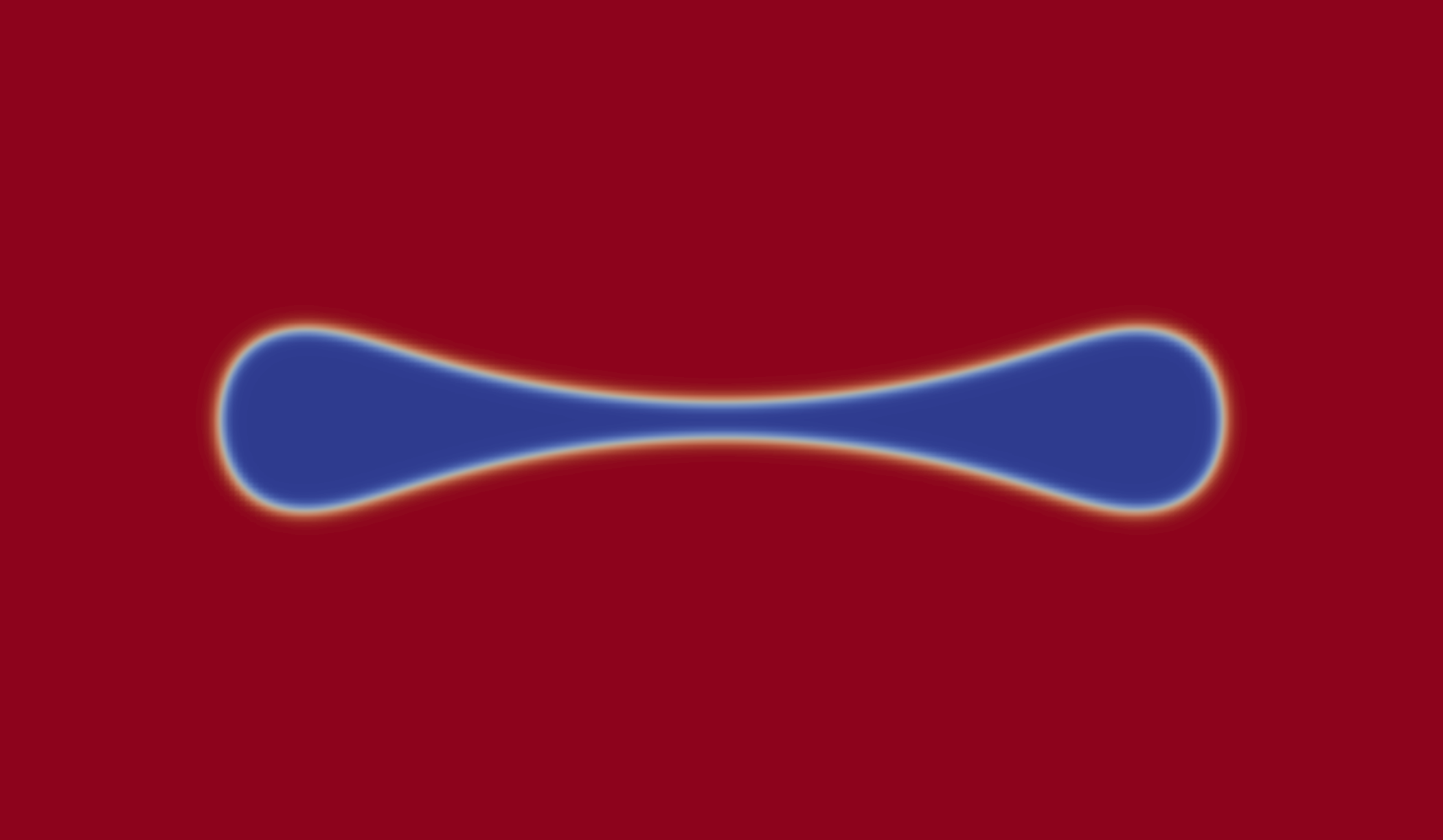}
\caption{{\color{black} Test 4. }The energy minimization method is used. The left graph is the initial condition, and the right graph is the evolution of the solution at around $t=0.00111$ {\color{black} with $\epsilon = 0.0033$}.}
\label{fig15}
\end{figure}

\begin{figure}[H]
\centering
\includegraphics[height=2.07in,width=2.3in]{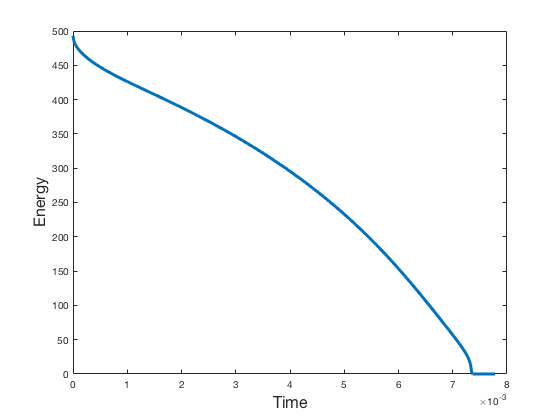}
\includegraphics[height=2.07in,width=2.3in]{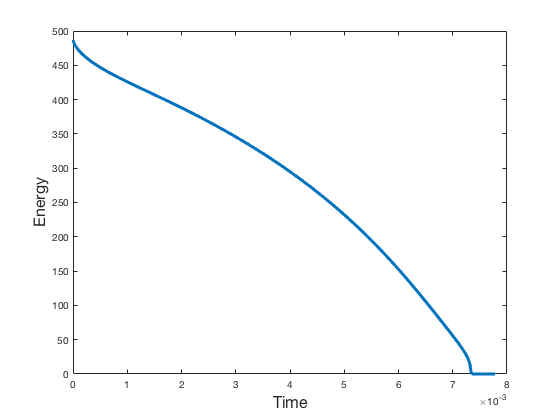}
\caption{{\color{black} Test 4. }The left hand side is the energy of the FIS, and the right hand side is the energy of the energy minimization method with $\epsilon = 0.0033$.}
\label{energy6}
\end{figure}

In Figure \ref{fig14}, the phase field method is used to compute the evolution of the solution, and we find that these two wedges merge, which is consistent with the result of the level set method. The spatial size is $h=0.0033$, the interaction length is $\epsilon=0.0033$, and the time step size is $k =1.11\times10^{-5}$. In Figure \ref{fig15}, {\color{black}using the same data for Figure \ref{fig14},} the energy minimization method is used to compute the evolution of the solution, and we find that these two wedges also merge, which is consistent with the result of the level set method. Figure \ref{energy6} indicates the energy change over time for both the phase field method and the energy minimization method. We observe that the energy decreases slowly until around time $7\times10^{-3}$ for both methods, then it decreases very fast until it reaches 0.  

\section{An Energy Penalized Minimization Algorithm.}\label{sec3}
In this section, an energy penalized minimization algorithm is proposed to compute the evolution of the mean curvature flow under topological changes. The penalty term has a very similar form as $E_n^{\rm AC}(u_h;u_h^{n-1})$, so it is very easy to implement this algorithm. The algorithm is defined by seeking $u_h^n$ such that
\begin{align}\label{our_scheme}
{u}_h^n = \underset{u_h\in V_h}{\mathrm{argmin}}\ \tilde E_n^{\rm AC}(u_h;u_h^{n-1}).
\end{align}
Here $\tilde E_n^{\rm AC}(u_h;u_h^{n-1})$ is defined by
\begin{equation}
\tilde E_n^{\rm AC}(u_h;u_h^{n-1}):=E_n^{\rm AC}(u_h;u_h^{n-1})+\delta(\frac{1}{\epsilon^2} F(u_h)-\frac{1}{\epsilon^2} F(u_h^{n-1})),
\end{equation}
 where $\delta>0$ is a penalty constant and $E_n^{\rm AC}$ is defined in \eqref{AC_energy}.
 
Next, we will prove the well-posedness and the energy stability of this energy penalized minimization algorithm.

\begin{theorem}\label{thm_add1}
Under the condition that $k\leq \frac{\epsilon^2}{1+\delta}$, the solution $u_h^n$ in \eqref{our_scheme} exists and is unique.
\end{theorem} 
\begin{proof}
By taking the second Fr$\acute{e}$chet derivative of $\tilde E_n^{\rm AC}(\cdot;u_h^{n-1})$, we
get for any $v_h \in V_h$,
\begin{align}\label{eq:FIS-AC-frechet}
(\tilde E_n^{\rm AC})''(u_h;u_h^{n-1})(v_h,v_h) =&
\frac{3(1+\delta)}{\epsilon^2}\int_{\Omega}u_h^2v_h^2dx \\
&\qquad+\int_{\Omega}(\frac{1}{k}-\frac{1+\delta}{\epsilon^2})v_h^2dx + \|\nabla
v_h\|_{L^2(\Omega)}^2.\notag
\end{align}
When $k\leq \frac{\epsilon^2}{1+\delta}$ and $v_h\neq 0$, $\tilde E(\cdot;u_h^{n-1})$ is strictly convex on $V_h$ since
\begin{equation*}
(\tilde E_n^{\rm AC})''(u_h;u_h^{n-1})(v_h,v_h)>0.
\end{equation*}
Then $u_h^n$ is the minimizer of a convex functional, so it exists and is unique.
\end{proof}

Instead of using the technique in Theorem \ref{lem20200421}, we use this minimization approach to show the energy stability property of $u_h^n$.
\begin{theorem}\label{thm_add2}
Under the condition that $k\leq \frac{\epsilon^2}{1+\delta}$, the following energy inequality holds
\begin{align*}
\tilde J_\epsilon^{\rm AC}(u_h^n)+ \frac{1}{{\color{black}2k}}\|u_h^n -u_h^{n-1}\|_{L^2(\Omega)}^2\le \tilde J_\epsilon^{\rm AC}(u_h^{n-1}),
\end{align*}
where $\tilde J_\epsilon^{\rm AC}(v_h):=J_\epsilon^{\rm AC}(v_h)+\frac{\delta}{\epsilon^2} F(v_h)$.
\end{theorem} 
\begin{proof}
Since $u_h^n$ is the global minimizer of the convex functional $\tilde E_n^{\rm AC}$, we have
\begin{align*}
&J_\epsilon^{\rm AC}(u_h^n) + \frac{1}{{\color{black}2k}}\|u_h^n -u_h^{n-1}\|_{L^2(\Omega)}^2+\delta(\frac{1}{\epsilon^2} F(u_h^n)-\frac{1}{\epsilon^2} F(u_h^{n-1}))\\
=& \tilde E_n^{\rm AC}(u_h^n;u_h^{n-1}) \\
\leq& \tilde E_n^{\rm AC}(u_h^{n-1};u_h^{n-1}) \\
=& J_\epsilon^{\rm AC}(u_h^{n-1}).
\end{align*}
Then the conclusion is obtained directly.
\end{proof}

In the following Test 5--Test 7, the energy penalized minimization algorithm is applied to implement the above benchmark problems as well as the problems with random initial conditions.\\

\begin{figure}[H]
\centering
\includegraphics[height=1.4in,width=2.3in]{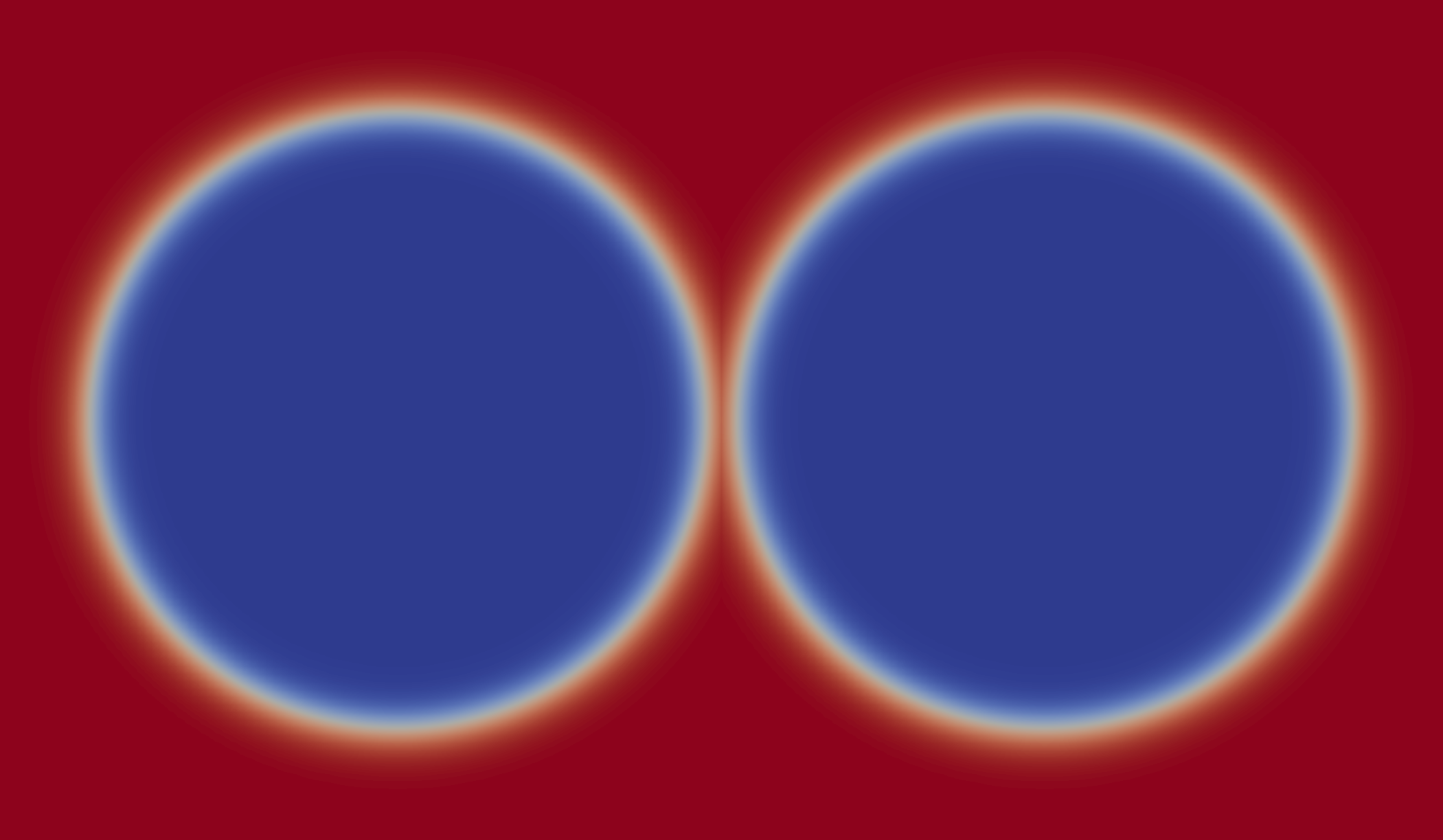}
\includegraphics[height=1.4in,width=2.3in]{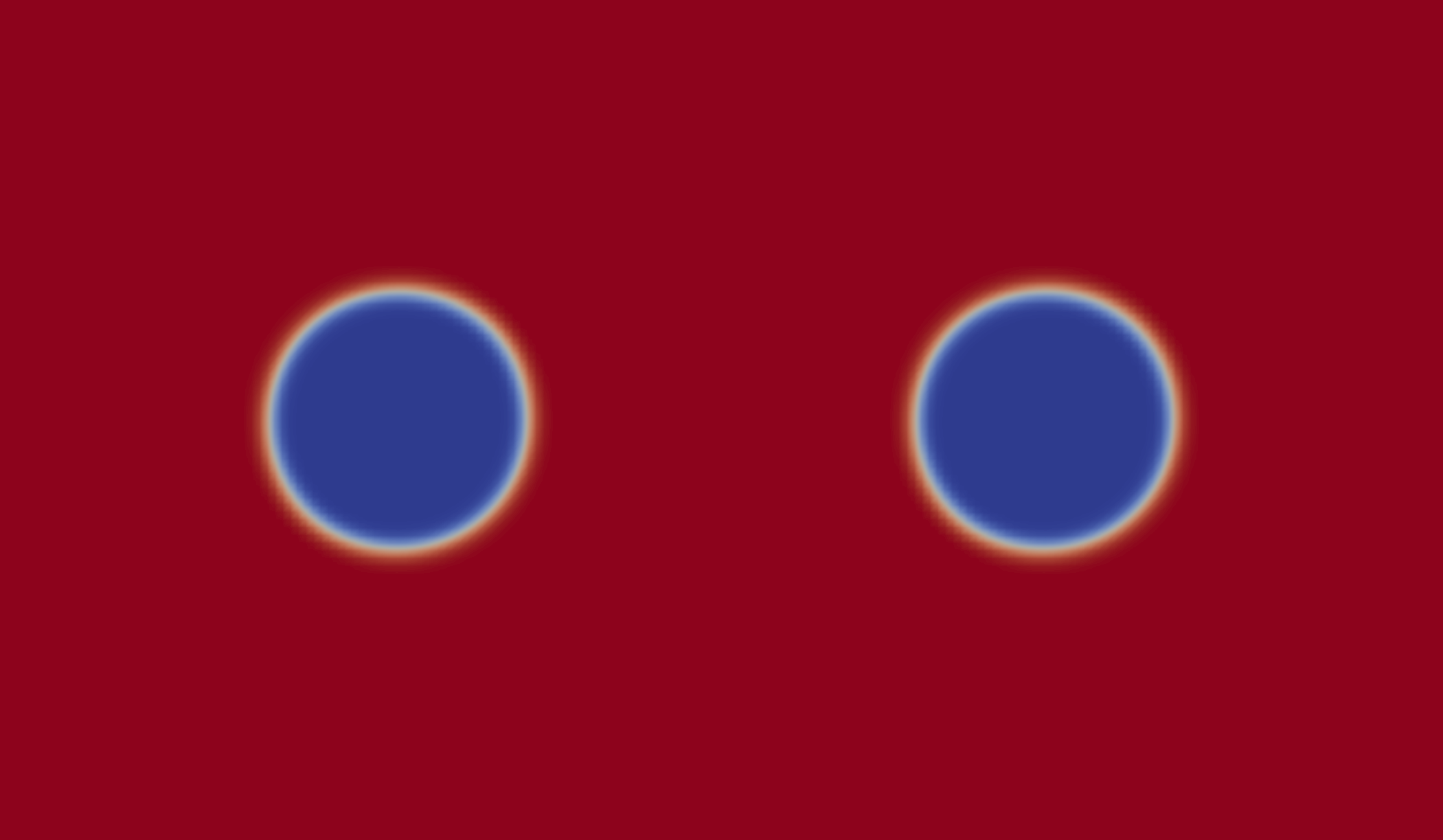}
\caption{{\color{black} Test 5. }The energy penalized minimization algorithm is used. The left graph is the initial condition, and the right graph is the evolution of the solution at $t=0.015$ {\color{black} with $\epsilon = 0.01$}.}
\label{fig10}
\end{figure}

{\bf Test 5.} %This should work well for two-circle case and random initial cases: Use new algorithm to implement Test 3 (including evolution and energy decay).
\label{test4} In this test, we compare the evolution of the solution when the initial condition is two circles with a small distance in Test 3 based on the energy penalized minimization algorithm. In Figure \ref{fig10}, the energy penalized minimization algorithm is used to compute the evolution of the solution, and we find that these two circles separate, which is consistent with the result of the level set method. The distance between these two circle is $d=0.02$, the spatial size is $h=0.005$, the interaction length is $\epsilon=0.01$, the penalized term is $\delta = 4$, and the time step size is $k = 1\times10^{-4}$. %The Figure \ref{energy4} indicates the energy change over time for the energy penalized minimization algorithm. We can observe that the energy decreases until it reaches 0 when these two circles disappear. 

%\begin{figure}[H]
%\centering
%\includegraphics[height=2.3in,width=3.4in]{test4_pic3.png}
%\caption{The energy plot of the energy penalized minimization algorithm.}
%\label{energy4}
%\end{figure}

\begin{figure}[H]
\centering
\includegraphics[height=1.4in,width=2.3in]{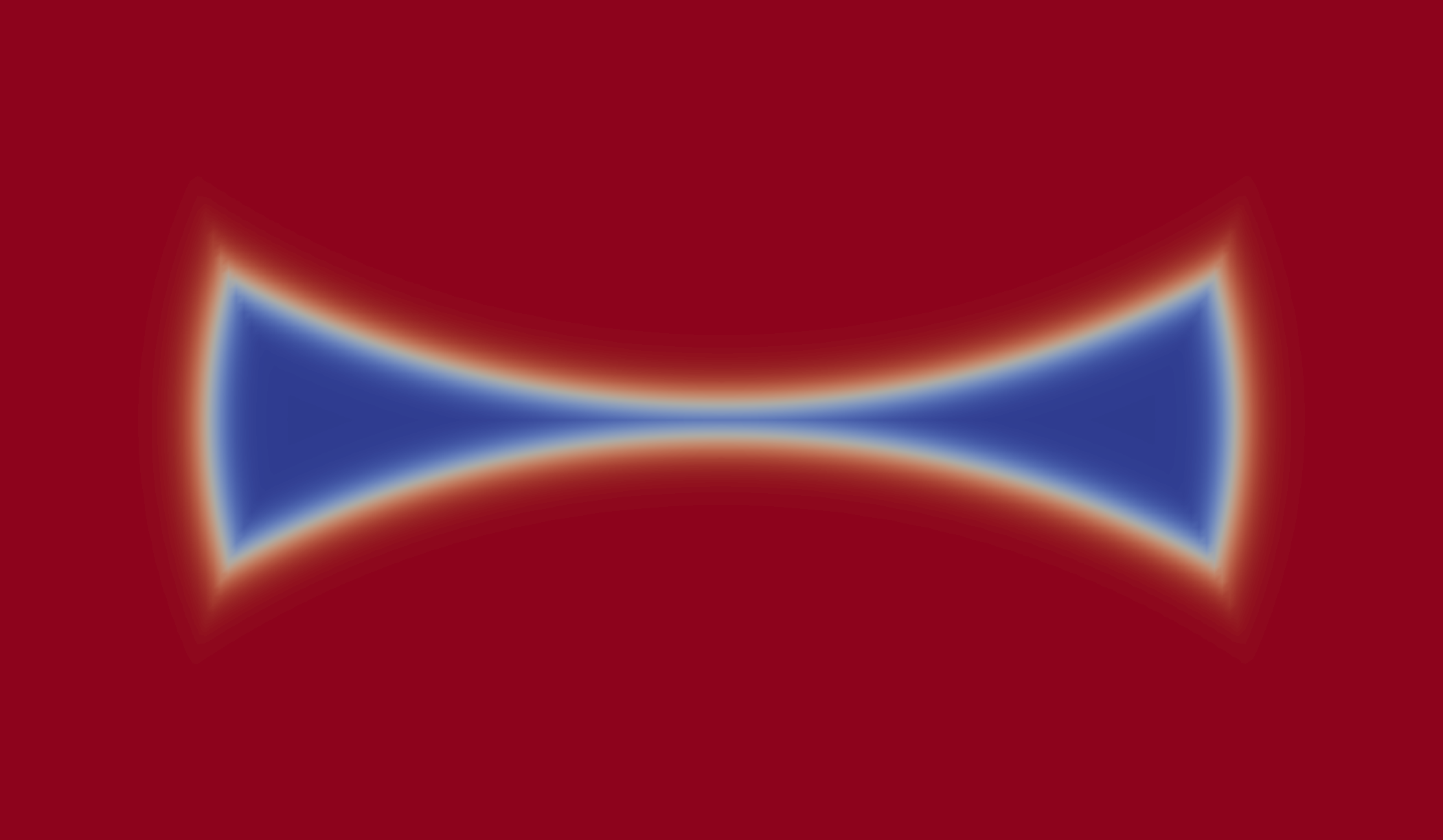}
\includegraphics[height=1.4in,width=2.3in]{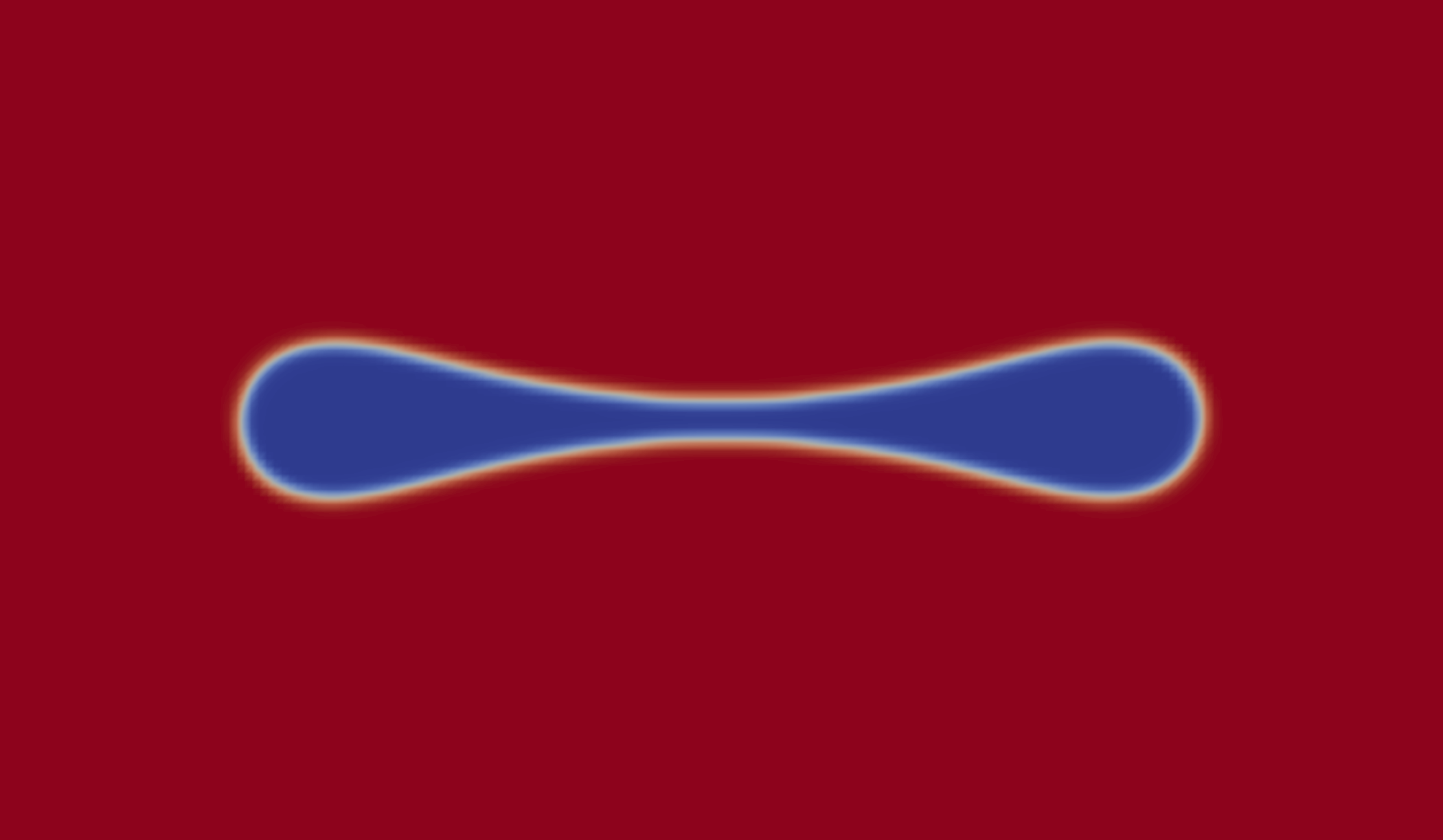}
\caption{{\color{black} Test 6. }The energy penalized minimization algorithm is used. The left graph is the initial condition, and the right graph is the evolution of the solution at $t=0.0015$ {\color{black} with $\epsilon = 0.01$}.}
\label{fig16}
\end{figure}

{\bf Test 6.}
In this test, we compare the evolution of the solution under the initial condition in Test 4 based on the energy penalized minimization algorithm. In Figure \ref{fig16}, the energy penalized minimization algorithm with penalized term $\delta$ to be 8 is used to compute the evolution of the solution, and we find that these two wedges merge, which is consistent with the result of the level set method. The spatial size is $h=0.005$, the interaction length is $\epsilon=0.01$, and the time step size is $k = 10^{-4}$. %The Figure \ref{energy7} indicates the energy change over time for the energy penalized minimization algorithm. We observe that the energy decreases slowly until around time $5\times10^{-3}$, then it decreases very fast until it reaches 0.  

\begin{figure}[H]
\centering
\subfloat[$t = 0$]{\includegraphics[height=1.9in,width=1.9in]{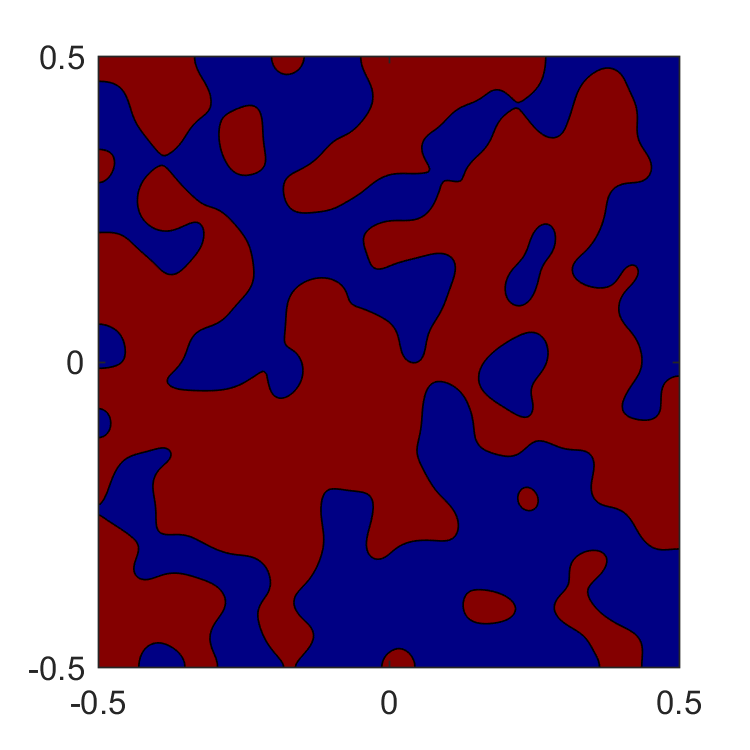}}
\subfloat[$t = 0.0005$]{\includegraphics[height=1.9in,width=1.9in]{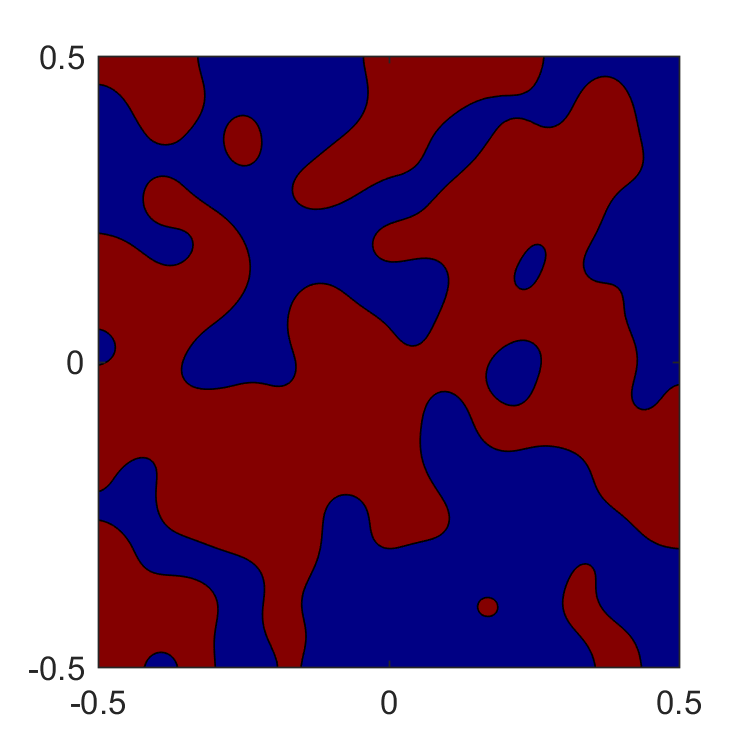}}

\subfloat[$t = 0.002$]{\includegraphics[height=1.9in,width=1.9in]{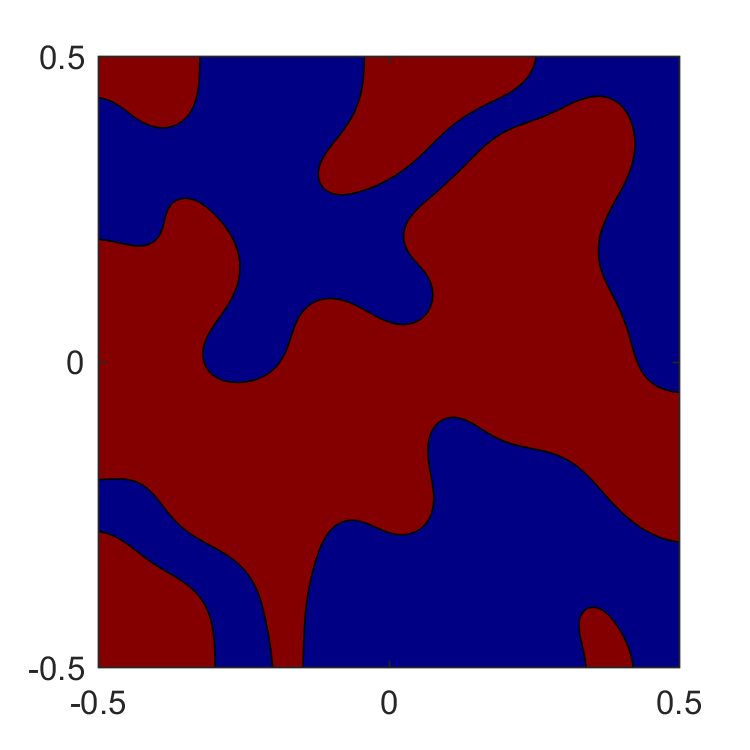}}
\subfloat[$t = 0.01$]{\includegraphics[height=1.9in,width=1.9in]{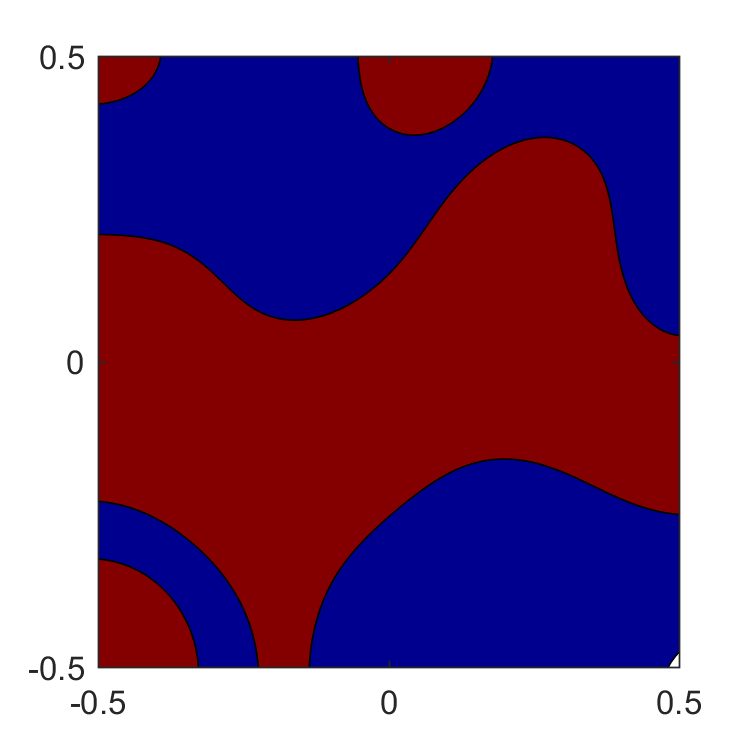}}
\caption{{\color{black} Test 7. }The level set method is used. The spatial size $h = 0.005$ and the time step size $k = 6.25 \times 10^{-6}$.}
\label{fig17}
\end{figure}

\begin{figure}[H]
\centering
\subfloat[$t = 0$]{\includegraphics[height=1.9in,width=1.9in]{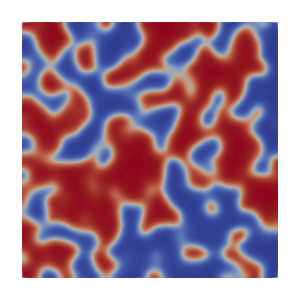}}
\subfloat[$t = 0.0005$]{\includegraphics[height=1.9in,width=1.9in]{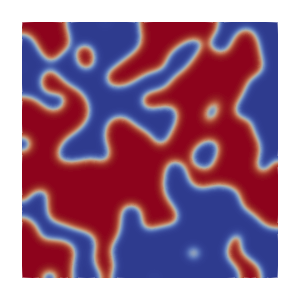}}

\subfloat[$t = 0.002$]{\includegraphics[height=1.9in,width=1.9in]{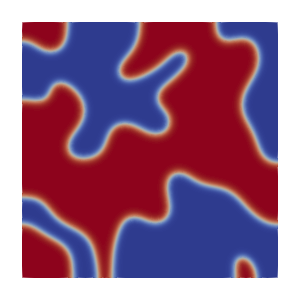}}
\subfloat[$t = 0.01$]{\includegraphics[height=1.9in,width=1.9in]{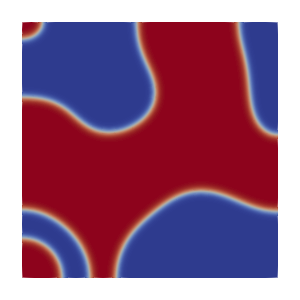}}
\caption{{\color{black} Test 7. }The energy minimization method is used. The spatial size $h = 0.005$ and the time step size $k = 5\times10^{-5}$ {\color{black} with $\epsilon = 0.01$}.}
\label{fig18}
\end{figure}

\begin{figure}[H]
\centering
\subfloat[$t = 0$]{\includegraphics[height=1.9in,width=1.9in]{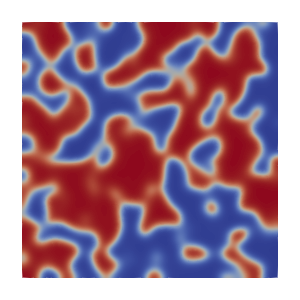}}
\subfloat[$t = 0.0005$]{\includegraphics[height=1.9in,width=1.9in]{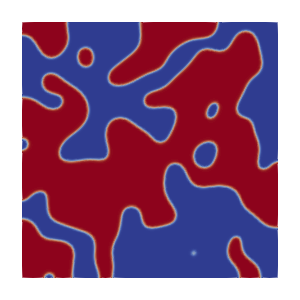}}

\subfloat[$t = 0.002$]{\includegraphics[height=1.9in,width=1.9in]{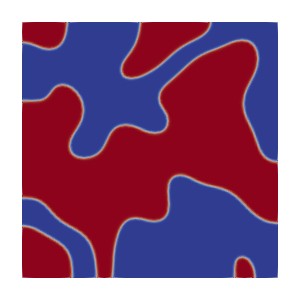}}
\subfloat[$t = 0.01$]{\includegraphics[height=1.9in,width=1.9in]{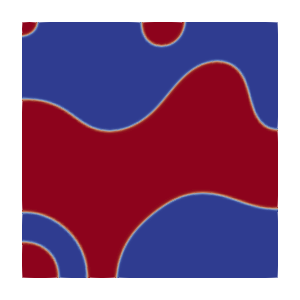}}
\caption{{\color{black} Test 7. }The energy penalized minimization algorithm is used. The spatial size $h = 0.005$ and the time step size $k = 5\times10^{-5}$ {\color{black} with $\epsilon = 0.01$}.}
\label{fig19}
\end{figure}

{\bf Test 7.}
In this test, we compare the evolutions of the solutions when the initial condition is random based on the level set method, the energy minimization method, and the energy penalized minimization algorithm. The domain is $\Omega = [-0.5, 0.5]^2$. We compute the evolution of the solution after 50 time steps when the initial condition is random, based on the phase field method with time step size $k = 10^{-5}$, and then use the resulted solution as the initial condition $u_0(x, y)$. We did many tests by choosing different random initial conditions, and we use one example to illustrate the advantages of the proposed energy penalized minimization algorithm in Figures \ref{fig17}--\ref{fig19}.
  
In Figure \ref{fig17}, the level set method is used to compute the evolution of the solution when the random initial condition is used, with the spatial size $h = 0.005$ and the time step size $k = 6.25 \times 10^{-6}$. In Figure \ref{fig18}, the energy minimization method is used. The spatial size is $h = 0.005$ and the time step size is $k = 5 \times 10^{-5}$. We observe that the evolution of the solution computed through the energy minimization is significantly different from the result of the level set method. In Figure \ref{fig19}, the energy penalized minimization algorithm is used. The spatial size is $h = 0.005$, the time step size is $k = 5 \times 10^{-5}$, and the penalized parameter $\delta = 8$. We observe that the evolution of the solution computed through the  energy penalized minimization algorithm is consistent with the pattern of the level set method. 

%Figure \ref{fig20} indicates the energy change over time for both the energy minimization method and the energy penalized minimization algorithm. We observe that, for both methods, the energy decreases very fast at the beginning; then the speed of decreasing slows down as time goes. However, the energy of the energy penalized minimization algorithm decreases faster than the energy minimization method.  
%
%\begin{figure}[H]
%\centering
%\includegraphics[height=2.5in,width=2.9in]{test7_pic13.png}
%\includegraphics[height=2.5in,width=2.9in]{test7_pic14.png}
%\caption{The left hand side is the energy of the energy minimization method, and the right hand side is the energy of the energy penalized minimization method.}
%\label{fig20}
%\end{figure}

%\begin{remark}
%From the figures of energy plots, we observe that the energy decreases fast when there are topological changes, i.e., Figure \ref{energy2} and Figure \ref{energy5}, and the energy decreases slowly when there are no topological changes, i.e., Figure \ref{energy1}, Figure \ref{energy3}, Figure \ref{energy6}, Figure \ref{energy4}, and Figure \ref{energy7}. However, the adaptive in-time technique, which uses energy as an error indicator to `prevent' the energy from decreasing fast, does not work for the initial conditions in this paper as mentioned in Test 3.
%\end{remark}

\section{A Multilevel Minimization Algorithm.}\label{sec4}
%%%%%%%%TO BE ADDED:
The coarse to fine multilevel minimization algorithm is defined by seeking, at each time step iteration, ${\color{black}{u}_{h_{m}}^{n}}$ such that
\begin{align*}
{u}_{{\color{black}h_1}}^{n-1} &= \underset{u_{{\color{black}h_1}}\in V_{{\color{black}h_1}}}{\mathrm{argmin}}\ E^{AC}_{n,{\color{black}\epsilon_1}}(u_{{\color{black}h_1}};u^{n-1}_{h_m}),\\
{u}_{{\color{black}h_2}}^{n-1} &= \underset{u_{{\color{black}h_1}}\in V_{{\color{black}h_2}}}{\mathrm{argmin}}\ E^{AC}_{n,{\color{black}\epsilon_2}}(u_{{\color{black}h_2}};u^{n-1}_{h_m}),\\
&\vdots\\
{\color{black}{u}_{h_{m}}^{n}} &= \underset{u_{{\color{black}h_m}}\in V_{{\color{black}h_m}}}{\mathrm{argmin}}\ E^{AC}_{n,{\color{black}\epsilon_m}}(u_{{\color{black}h_m}};u^{n-1}_{h_m}),
\end{align*}
where $h_1>h_2>\cdots>h_m$, $\epsilon_1>\epsilon_2>\cdots>\epsilon_m$, $V_{h_i}$ is the discrete space associated to $h_i$, and $E^{AC}_{n,{\color{black}\epsilon_i}}$ is the energy defined in \ref{AC_energy} using $\epsilon_i$ instead of $\epsilon$. This means that for the first few steps , when $\epsilon_i$ and $h_i$ are large, our energy minimization should converge quickly. Then, when $\epsilon_i$ and $h_i$ are smaller, the minimization algorithm is able to converge due to the fact that the initial guess  which is the solution from the previous mesh is close enough to the physical solution. Note that this method requires to project the solution obtained from the previous step to the new mesh; as for the minimization algorithm we use the limited memory BFGS method, see \cite{nocedal1980}.
%%%%%%%%%%

One of the issues of using the traditional energy minimization method when the energy is not convex is that the solution is dependent on the initial guess. The method we propose to solve this issue is, at a fix time step, first minimizing the energy with a larger $\epsilon$ in order to start with a convex energy functional, and then decreasing $\epsilon$ in the next steps. Regarding the fully implicit scheme for the Allen-Cahn equation, equation \eqref{FIS-AC}, we recall the following two important conditions we need:

\begin{enumerate}
 \item Convexity condition: ${\color{black}k} \le \epsilon ^2$
 \item Mesh size relation to $\epsilon$: $h^{-1}=\mathcal
O(\epsilon^{-1})$
\end{enumerate}

In the following Test 8 and Test 9, the traditional energy minimization algorithm and the coarse to fine multilevel minimization algorithm are compared for the cases when there are no topological changes (Test 8) and there are topological changes (Test 9).\\

{\bf Test 8.} To test this method we pick the following: $h=10^{-3}$, $\epsilon=0.005$, and
\begin{equation}
 u_0(x,y) = \tanh\bigl(\frac{d_0(x,y)}{\sqrt{2}\epsilon}\bigr).
\end{equation}
For our initial guess we chose $u = 1-u_0(x,y)$.
%\begin{equation}
% u = 1-u_0(x,y).
%\end{equation}
We remark that when $\epsilon$ is large, by following the condition $h =\mathcal O(\epsilon)$, we can pick a larger mesh size. This leads up to the following multilevel algorithms:

\begin{figure}[h]
 \centering
 \begin{tabular}{ccccccc} 
  \includegraphics[scale=0.079]{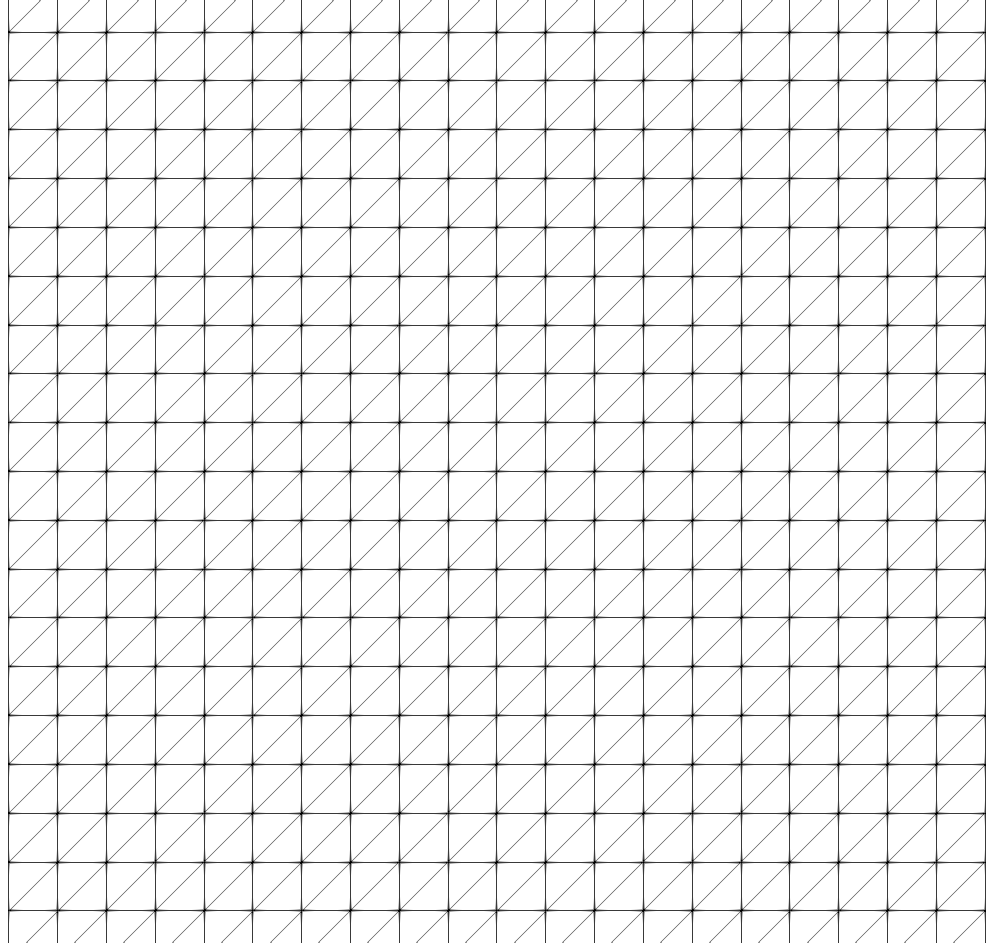} &  &\includegraphics[scale=0.079]{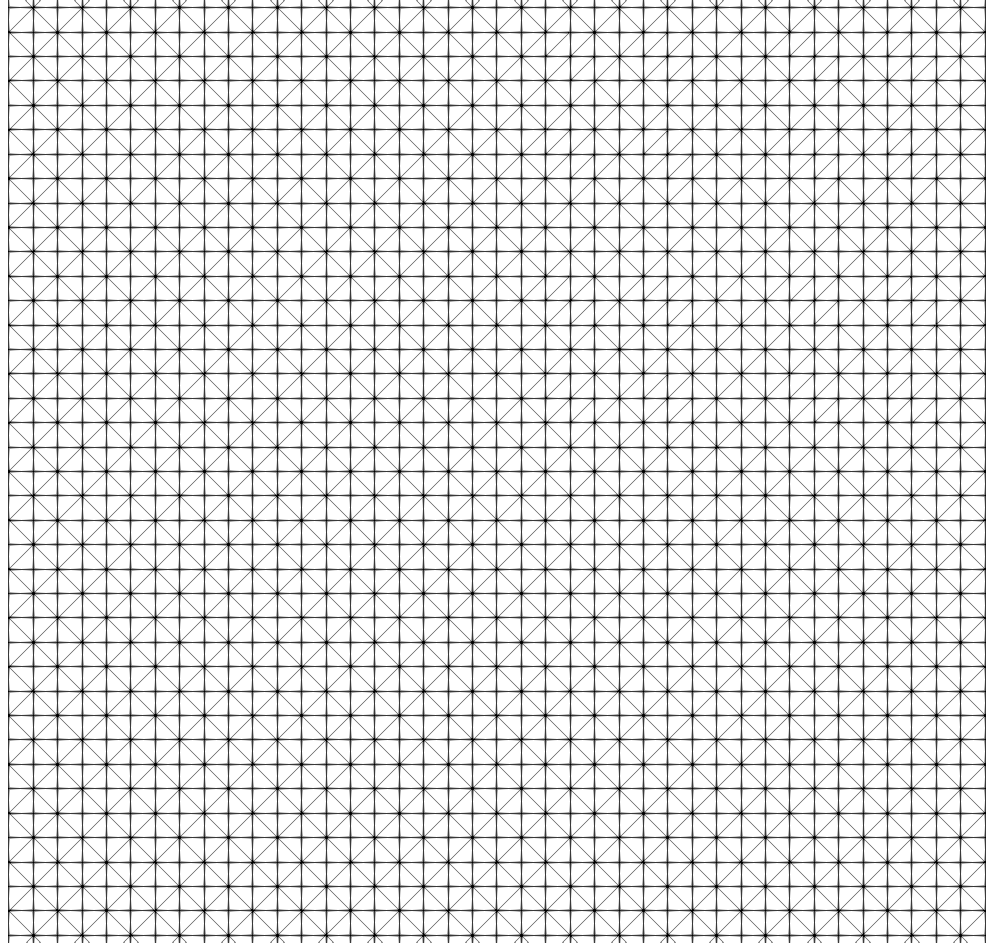} &  & \includegraphics[scale=0.079]{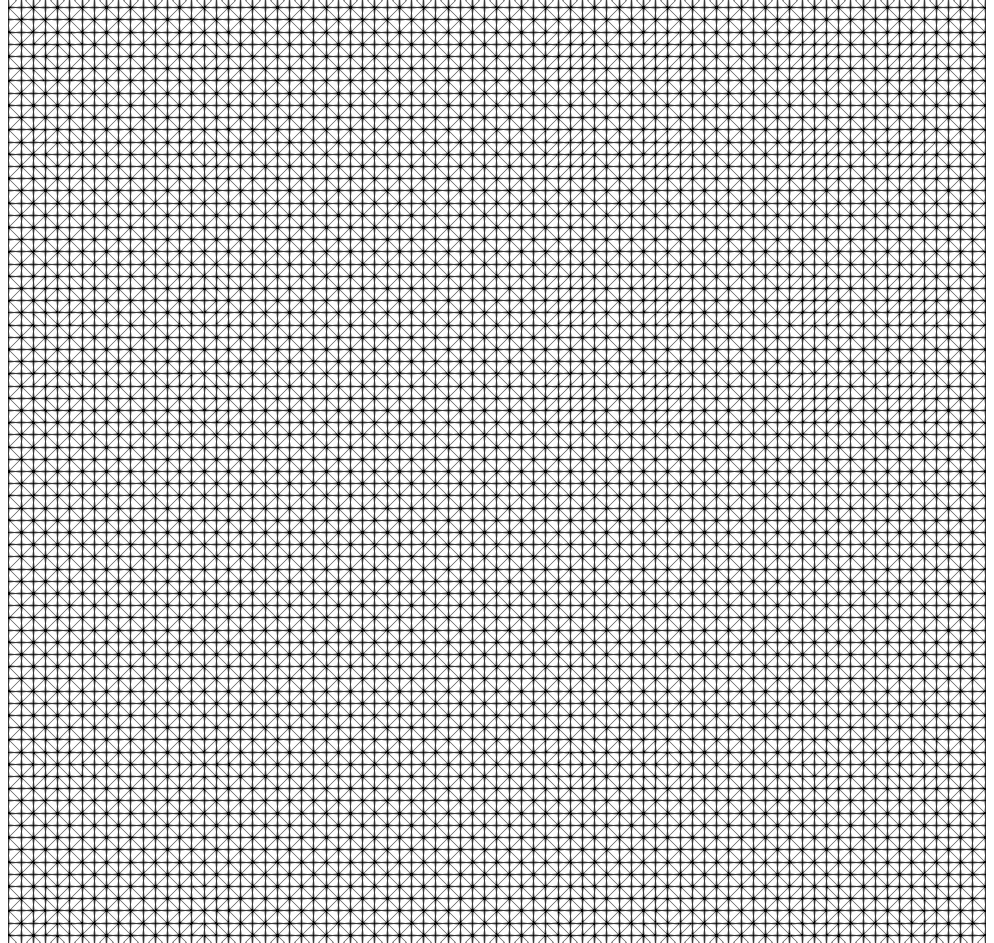} \\
  $h=0.12$, $\epsilon=0.1$ & $\longrightarrow$ & $h=0.07$, $\epsilon=0.076$ & $\longrightarrow$ & $h=0.035$, $\epsilon=0.028$ & $\cdots$ &%$\longrightarrow$ & $\cdots$
 \end{tabular}
\end{figure}

We chose 5 different meshes with the following $h$ and $\epsilon$ in Table \ref{tab1}.
\begin{table}[!htbp]
\centering
 \begin{tabular}{l|l|l}
 Mesh & $h$ & $\epsilon$ \\
 \hline\hline
 1 & 0.12 & 0.1 \\
 2 & 0.07 & 0.076 \\
 3 & 0.035 & 0.28\\
 4 & 0.018 & 0.29\\
 5 & 0.009 & 0.005
 \end{tabular}
 \caption{Meshes on different levels.}\label{tab1}
\end{table}

\begin{figure}[H]
\centering
\includegraphics[scale=0.31]{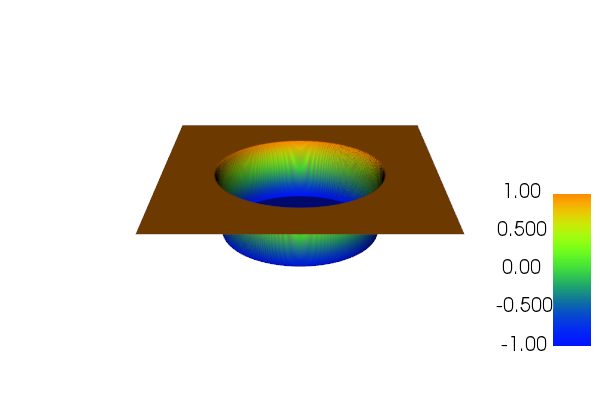}
\includegraphics[scale=0.31]{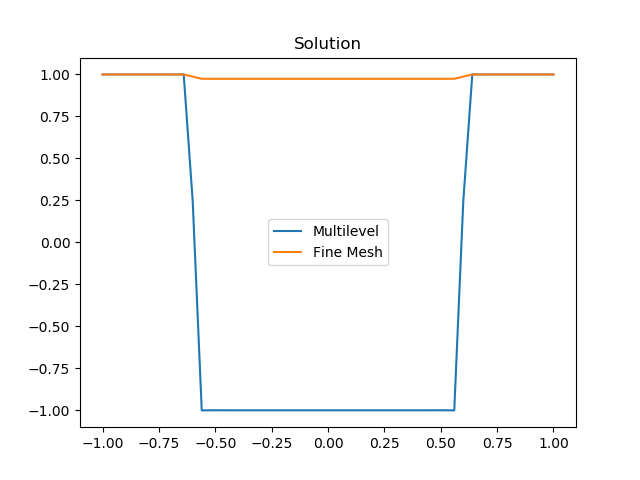}
\caption{{\color{black} Test 8. }Left: Initial value of $u$ for the Allen-Cahn equation. Right: Cross-sectional solutions to the Allen-Cahn equation at $y=0$ and $t=0.001$ of the regular and multilevel methods.}
\label{fig4.2a}
\end{figure}

Then, after doing one step in time we compare the multilevel solution and the solution computed by using the fine mesh only. Figure \ref{fig4.2a} (right) displays the cross-sectional solutions at $y = 0$ at $t=0.001$. We know, from \cite{xu2019stability}, that the correct solution is a circle of depth $-1$, see Figure \ref{fig4.2a} (left). Figure \ref{fig4.2a} (right) shows that the solution computed using the multilevel algorithm is the correct one. Notice the energy functional is convex for the coarse meshes of the multilevel minimization algorithm, but it is not convex for the fine meshes. We observe that the proposed multilevel minimization algorithm is useful to handle the issue of the initial guess when the energy functional is not convex.

{\bf Test 9.}
In this test, we compare the evolutions of the solutions based on the multilevel minimization algorithm and the traditional energy minimization algorithm \eqref{tra_energy} under the initial condition of two circles with a small distance. The distance between two circles is $d=0.02$, the spatial size is $h = 0.002$, the time step is $k = 1 \times 10^{-4}$, and the interaction length is $\epsilon=0.002$. The initial guess we choose is: $ u = 1-u_0(x,y).$
%\begin{equation}
% u = 1-u_0(x,y).
%\end{equation} 

\begin{figure}%[h]
 \centering
 \includegraphics[height=1.4in,width=2.3in]{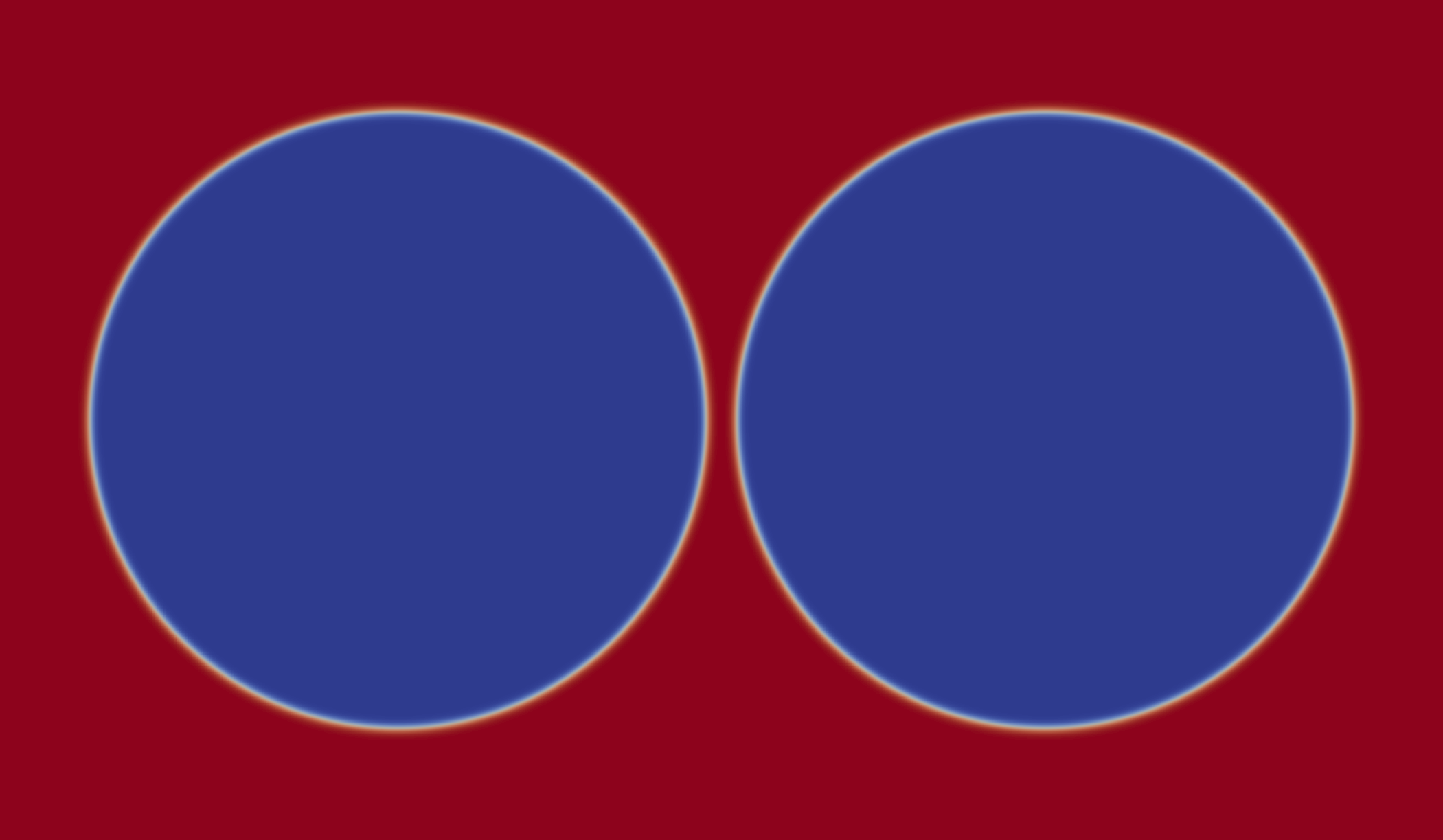}
\includegraphics[height=1.4in,width=2.3in]{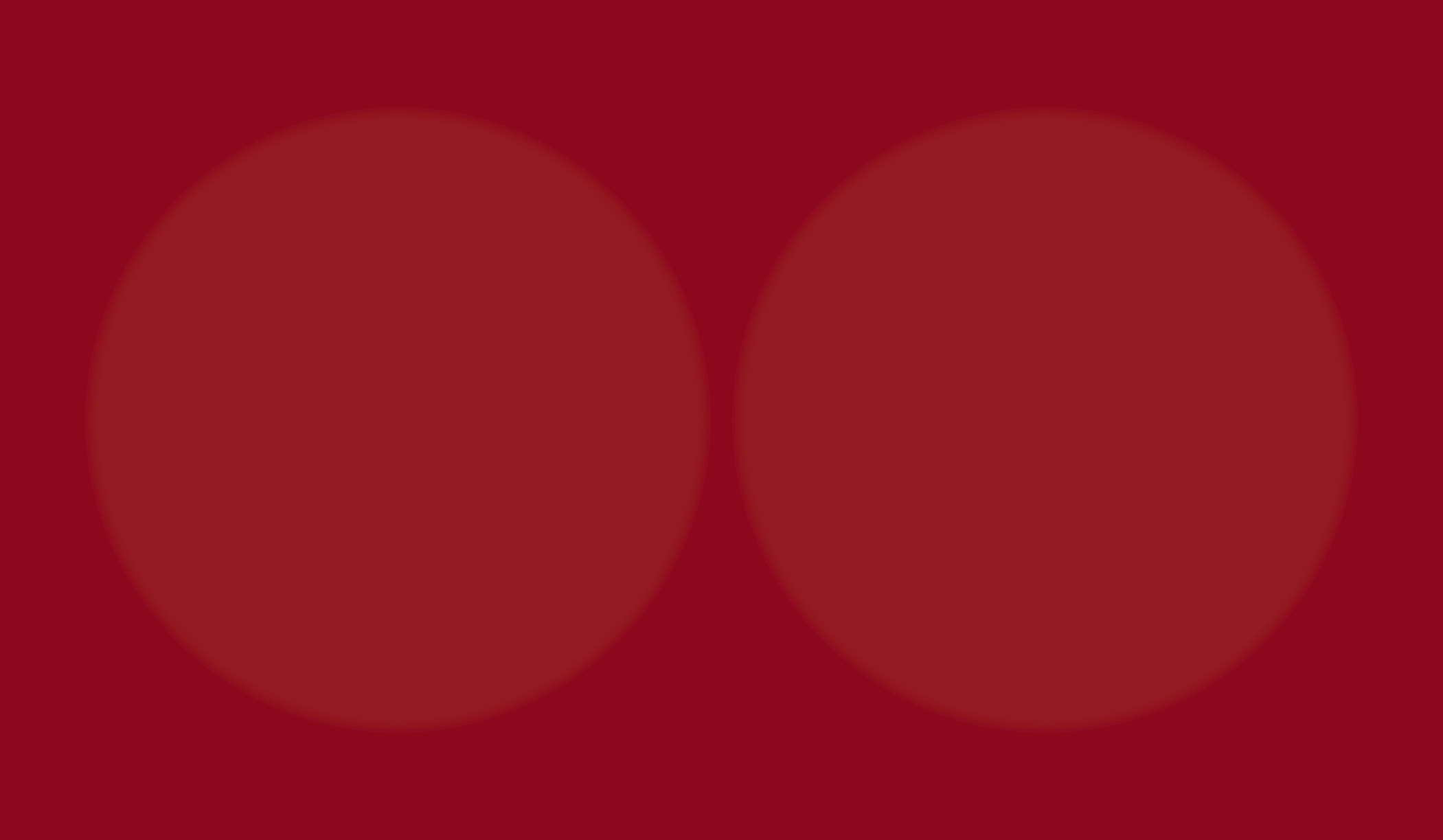}
\caption{{\color{black} Test 9. }The traditional energy minimization algorithm is used. The graph on the left is the initial condition, and the graph on the right is the evolution of the solution at $t = 0.0001$ with $\epsilon = 0.002$.}
 \label{fig4.2c}
\end{figure}

Figure \ref{fig4.2c} shows the evolution of the solution based on the traditional energy minimization algorithm. We observe that two circles become two round-shape shadows at $t = 0.0001$, which is not correct.

\begin{figure}[H]
 \centering
 \includegraphics[height=1.4in,width=2.3in]{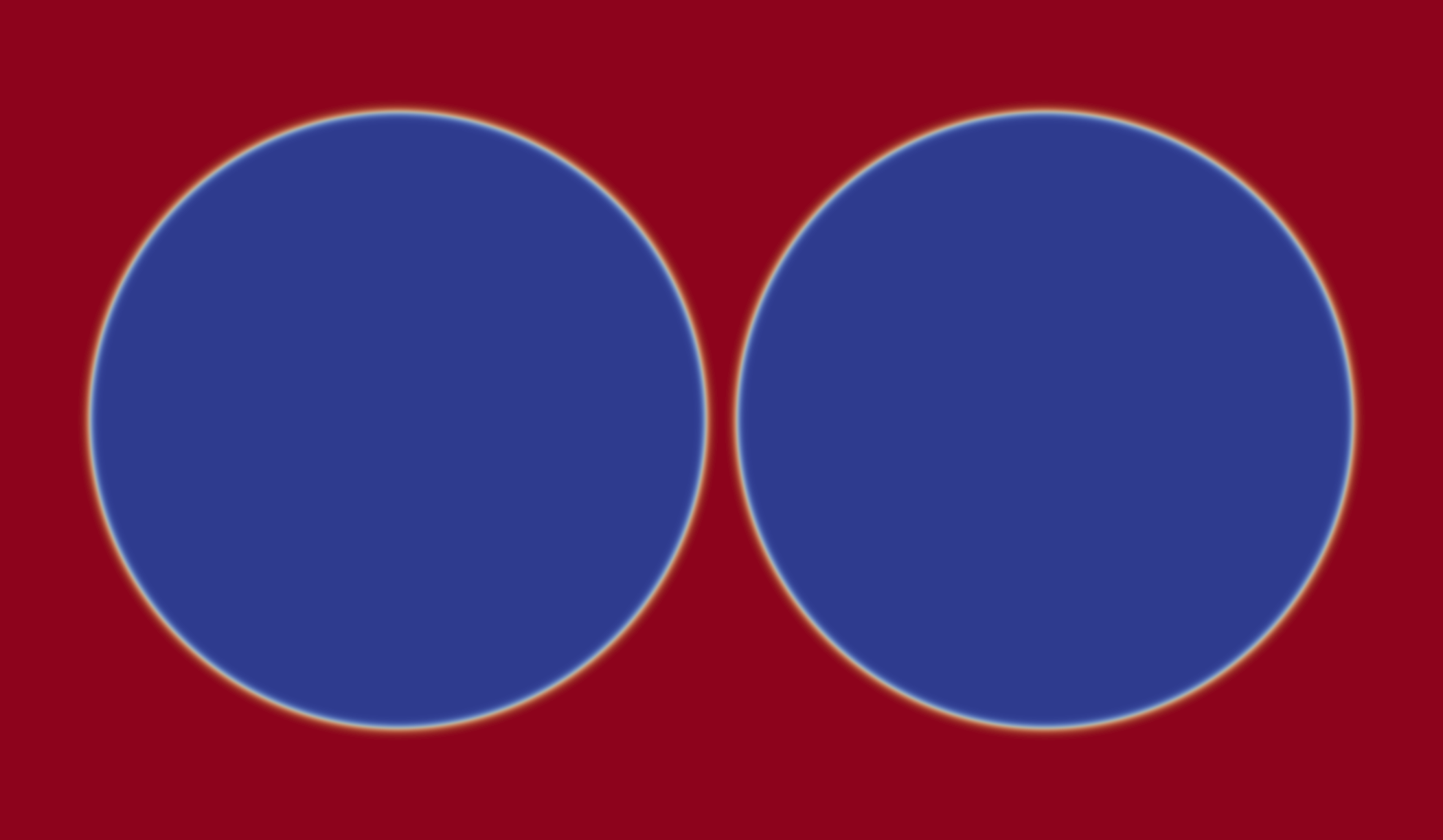}
\includegraphics[height=1.4in,width=2.3in]{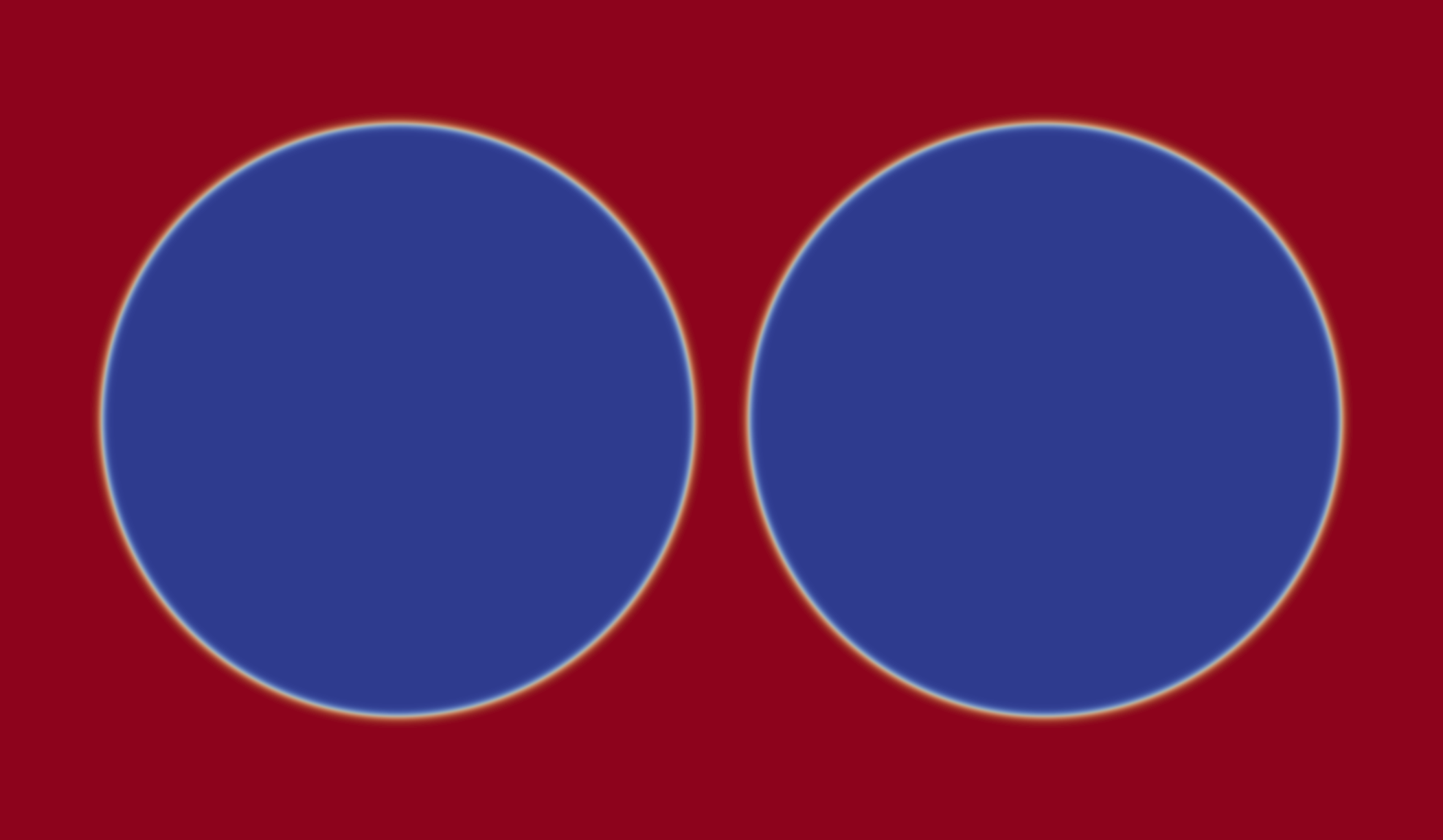}
\caption{{\color{black} Test 9. }The multilevel minimization algorithm algorithm is used. The graph on the left is the initial condition, and the graph on the right is the evolution of the solution at $t = 0.0015$.}
 \label{fig4.2d}
\end{figure}

Figure \ref{fig4.2d} shows the evolution of the solution based on the coarse to fine multilevel minimization algorithm. We observe that these two circles separate, which is consistent with the level set method. The distance between these two circles is $d=0.02$ and the time step size is $k = 1\times 10^{-4}$. Specifically, for the multilevel method, we pick $\epsilon$ and the spatial size $h$ in Table \ref{tab2} below: \\
\begin{table}[!htbp]
\centering
 \begin{tabular}{l|l|l}
 Mesh & $h$ & $\epsilon$ \\
 \hline\hline
 1 & 0.002 & 0.0047 \\
 2 & 0.001 & 0.00335 \\
 3 & 0.0005 & 0.002
 \end{tabular}
  \caption{Meshes on different levels}\label{tab2}
\end{table}

This test indicates that the traditional energy minimization algorithm may not be able to carry out the evolution of the solution correctly when the initial guess is not good and when there are topological changes. But the multilevel minimization algorithm may overcome the problem and produce the correct results.

\section{Conclusion.}\label{sec5}
In this paper, we propose some algorithms to compute the mean curvature flow under topological changes based on the phase field methodology, and the objective is to validate these algorithms for the random initial conditions. To achieve this objective, some benchmark problems are constructed first, and it is shown that the evolutions of the solutions are very sensitive to the interaction length $\epsilon$ when topological changes happen. The energy penalized minimization algorithm is proposed to accurately solve these benchmark problems and the problems with random initial conditions in Section \ref{sec3}. Besides, the issue of traditional energy minimization algorithm arises from the existence of multiple minimum points due to the non-convexity property of the energy functional. This issue can be removed by using a multilevel algorithm and starting from a convex problem in Section \ref{sec4}. With the multilevel algorithm, the solution converges to the global minimum. Furthermore, in Test 9, we show that when topological changes happen, the traditional energy minimization algorithm also encounters the convexity problem and fails to produce a correct solution. Fortunately, the multilevel algorithm is not only able to remove the convexity problem, but also able to carry out the evolution of the solution accurately. 

Lastly, we are also thinking about how to improve even further our energy minimization algorithms. We could think about using the multilevel methods inside the nonlinear solver. Instead of doing all the minimization steps on coarse meshes, we could first do some steps on the coarser meshes and then project the solution on the finer meshes. This could in theory reduce the time needed by our nonlinear solver.

\end{document}